\documentclass[10pt]{article}
\usepackage{combgames}
\usepackage{hyperref}
\usepackage{enumerate}
\usepackage{geometry}
\usepackage{amsmath,amsfonts,amsthm,amssymb}

\newtheorem{theorem}{Theorem}[section]
\newtheorem{proposition}[theorem]{Proposition}
\newtheorem{definition}[theorem]{Definition}
\newtheorem{lemma}[theorem]{Lemma}
\newtheorem{conjecture}[theorem]{Conjecture}
\newtheorem{corollary}[theorem]{Corollary}

\newtheorem{example}[theorem]{Example}
\newtheorem{question}[theorem]{Question}
\newtheorem{problem}[theorem]{Problem}

\newcommand{\cgname}[1]{{\sc #1}}
\newcommand{\cgbgame}[2]{\combgame{\{#1\mid#2\}}}

\newcommand{\qm}{s}
\newcommand{\qp}{q}
\newcommand{\R}{\mathcal{R}}
\newcommand{\N}{\mathbb{N}}
\newcommand{\Nz}{\mathbb{N}_0}
\newcommand{\n}{\mathcal{Q}}

\newcommand{\WPS}{\cgname{GoldenNugget} }
\newcommand{\wps}{\text{GN}}

\newcommand{\nug}{\text{heap\ }}

\newcommand{\heap}[1]{\ensuremath{{\rm  val}({#1})}}
\newcommand{\rcf}[1]{\ensuremath{{\rm rcf}({#1})}}
\newcommand{\cf}[1]{\ensuremath{{\rm  val}({#1})}}

\renewcommand{\hat}{\widehat}
\renewcommand{\ge}{\geqslant}
\renewcommand{\geq}{\geqslant}
\renewcommand{\le}{\leqslant}
\begin{document}

\title{Finding Golden Nuggets by Reduction}
\author{Urban Larsson\footnote{Department of Mathematics and Statistics, Dalhousie University, PO Box 15000, Halifax, Nova Scotia, Canada, B3H 4R2, urban031@gmail.com; supported by the Killam Trust} \and Neil A. McKay\footnote{Department of Mathematics and Statistics, Dalhousie University, PO Box 15000, Halifax, Nova Scotia, Canada, B3H 4R2, nmckay@chebucto.ns.ca}, \and Richard J. Nowakowski\footnote{Department of Mathematics and Statistics, Dalhousie University, PO Box 15000, Halifax, Nova Scotia, Canada, B3H 4R2, rjn@mathstat.dal.ca}, \and Angela A. Siegel\footnote{Robert Gordon University, Garthdee House, Garthdee Road, Aberdeen, AB10 7QB, Scotland, UK, a.a.siegel@rgu.ac.uk}}
\maketitle

\begin{abstract}
We introduce a class of normal play partizan games, called Complementary Subtraction. Let $A$ denote your favorite set 
of positive integers. This is Left's  subtraction set,  whereas Right subtracts  numbers not in $A$. The Golden Nugget Subtraction Game has the $A$ and $B$ sequences, from Wythoff's game, as the two complementary subtraction sets. As a function of the heap size, the maximum size of the canonical forms grows quickly. However, the value of the heap is either a number or, in reduced canonical form, a switch. We find the  switches by using properties of the Fibonacci word and standard 
Fibonacci representations of integers. Moreover, these switches are invariant under shifts by certain Fibonacci 
numbers. The values that are numbers, however, are distinct, and we find a polynomial time bit characterization for them, 
via the ternary Fibonacci representation.
\end{abstract}

\section{Introduction}
The game of \textsc{WythoffNim} \cite{W} is one of the earliest combinatorial games to be analyzed. A recent trend in the theory of \emph{impartial subtraction games}, starting with \cite{LHF, UL}, is to study new games defined via the set of P-positions (previous player winning positions) in the old game as subtraction sets in the new game. We study a \emph{partizan subtraction game} related to the P-positions of \textsc{WythoffNim}. To our knowledge, this is the first such study for partizan games.

\textsc{WythoffNim} is \emph{impartial}, meaning that the set of options does not depend on whose turn it is. The analysis of \textsc{WythoffNim} involves two \emph{complementary} sequences (each positive integer is in precisely one of the sequences). Let $\phi = \frac{\sqrt{5}+1}{2}$ be the \emph{golden ratio}. Wythoff's sequences are $A(n)=\lfloor  n\phi\rfloor $ and $B(n)=\lfloor  n\phi^2\rfloor $, for $n\in \N$, where $\N$ denotes the positive integers. They occur together as $\{(A(n), B(n)), (B(n), A(n))\mid n\in \Nz \}$ which is the set of P-positions of \textsc{WythoffNim}, where $\Nz =\N\cup \{0\}$ and $A(0)= B(0)=0$. Wythoff \cite{W}, needing a non-calculator means of generating the sequences, uses the pair of recurrences for integers $n\ge 0$: $A(n) = \mbox{mex}\;\{A(i),B(i)\mid i < n\},$ where $\mbox{mex}\,X$ denotes the least nonnegative number not in $X\subset \Nz$ (a proper subset) and 
\begin{equation}
B(n) = A(n) + n.\label{bann}
\end{equation}

\begin{table}[ht]
\begin{center}
\begin{tabular}{|l||c|c|c|c|c|c|c|c|c|c|c|c|c|c|c|}\hline
$n$      & 0& 1& 2& 3& 4& 5& 6& 7& 8& 9&10&11&12&13&14\\ \hline\hline
$A(n)$    & 0& 1& 3& 4& 6& 8& 9&11&12&14&16&17&19&21&22\\ \hline
$B(n)$    & 0& 2& 5& 7&10&13&15&18&20&23&26&28&31&34&36\\ \hline
\end{tabular}\caption{The first few terms of the $A$ and $B$ sequences.}\label{tab:sequences}
\end{center}
\end{table}

\textsc{WythoffNim} and also \textsc{Nim} \cite{B} belong to the class of \emph{subtraction games}, of which the most widely known 
 game is probably the children's game of \textsc{TwentyOne}, where two players alternate to subtract either one or two from the given 
 number and starting with 21; a player who cannot move loses\footnote{Also seen in the television shows, \textit{Sesame Street} 
 and \textit{Survivor}.}. Guy \cite{GSm} uses the Sprague-Grundy theory for impartial games to analyze many more subtraction 
 (and also splitting) games. In particular, it is easy to prove that the value-sequence 
 of any (one heap) impartial subtraction game on a finite subtraction set,  is eventually periodic. There are many examples of related regularities known in games: e.g., value sequences 
 that are periodic, arithmetic-periodic, and split-periodic \cite{ANW}, and the ruler regularity sequences \cite{GN}. 

Wythoff's sequences are among the first  known that are aperiodic
and regular. Exploring their properties, e.g \cite{S1, K}, and generalizing the golden ratio to other irrationals
 have spawned an industry. For example, Beatty sequences \cite{Be} and, of special interest in this paper, Sturmian words and the Fibonacci morphism \cite{L}, explained in Section \ref{sec:fibW}.
 Having two sequences suggests a \emph{partizan} subtraction game; the players get different \emph{subtraction sets}. 

\begin{definition}[\WPS\!\!] Consider a heap with a nonnegative number of counters. A Left option is to remove any number of tokens, provided it belongs to \textsc{WythoffNim}'s $A$ sequence. A Right option is to remove any number, provided it belongs to the $B$ sequence. The first player who cannot move, loses.
\end{definition}

For example, with 5 counters, Left can change the heap to one of size 1 (subtract 4),
2 (subtract 3), or 4 (subtract 1) whereas Right can change it to 0 or 3. For brevity, we will abbreviate \WPS to \wps.

Fraenkel and Kotzig
\cite{FK} show that for partizan subtraction games, with finite subtraction sets,
 the outcome sequence is eventually periodic. Here, it is easy to see that this is not the case. If the heap size is in $A$ then Left wins
 and if it is in $B$ then it is a first player win and the position need not last more than two moves. (See Theorem \ref{thm:outcomes}.) 
 In general, then, the one heap game is not very interesting. 
So let's play with more than one heap. 
 The heaps are composed of either blue or red counters. In `blue' heaps Left can remove any number in the $A$ sequence but in
 `red' heaps Right has this privilege.
 For example, consider the position with a blue heap of 9 and a red heap of 5 which we will denote $9_b, 5_r$.
 Left can change the blue heap to any of $1_b$,   $3_b$,     $5_b$,  $6_b$, and $8_b$, and the red heaps remains the same; or  she can leave the blue heap intact and change the red heap to either of $0$ and $3_r$.  (For those game-playing readers, it is clear that Left moves to $5_b$ for an easy win.)
  
\begin{question}\label{qu:example}  Consider the positions (i) with blue heaps of 3 and 20 and a red heap of 18; (ii) a blue heap of 20 and a red heap of 17.
Who wins these positions? See the last section for the answers.
\end{question}

To determine the winner, the outcome sequences give little information about the outcome when there is more than one heap.
However, values help us understand the game; the value gives the `number of moves' advantage, positive for Left 
and negative for Right.

We write $\cf{h}$ when referring to the \emph{canonical form} (value) of the \nug of size $h$ \cite{C}. 
We write $\rcf{h}$ for the \emph{reduced canonical form} of the \nug of size $h$. The reduced canonical form of a position, $G$, is the simplest position infinitesimally close to $G$; thus allowing us to give an analysis by ignoring infinitesimals \cite{GSi}. We give a formal introduction to reduced canonical form and other relevant game theory in Section \ref{sec:game theory}.

 \begin{table}[ht!]
\begin{centering}
\begin{tabular}{|c||l|l|}\hline
$h$&value & $\rcf{h}$\\\hline\hline
1&1&1\\ [2pt]
2&$\{1|0\}$&$\{1|0\}$\\ [2pt]
3&$\frac{1}{2}$&$\frac{1}{2}$\\ [2pt]
4&$\{1||1|0\}$&1\\ [2pt]
5&$\{1,\{1|0\}|0\}$&$\{1|0\}$\\ [2pt]
6&$\frac{3}{4}$&$\frac{3}{4}$\\[2pt]
7&$\{1||1|0|||0,\{1|0\}\}$&$\{1|0\}$\\[2pt]
8&$\{1|\frac{1}{2}\}$&$\{1|\frac{1}{2}\}$\\[2pt]
9&$\{1|\{1|0\},\{1||1|0\}\}$&1\\[2pt]
10&$\{1,\{1|0\}|0,\{1,\{1|0\}|0\}\}$&$\{1|0\}$\\[2pt]
11&$\frac{5}{8}$&$\frac{5}{8}$\\ [2pt]
12&$\{1||1|0||||1||1|0|||0,\{1|0\}\}$&$1$\\ [2pt]
13&$\{1,\{1,\{1|0\}|0\}|0\}$&$\{1|0\}$\\ [2pt]
14&$\frac{7}{8}$&$\frac{7}{8}$\\ [2pt]
15&$\{\{1||1|0\},\{1||1|0|||0,\{1|0\}\}|0,\{1|0\}\}$&$\{1|0\}$\\ [2pt]
16&$\{1,\{1|\frac{1}{2}\}|\frac{1}{2}\}$&$\{1|\frac{1}{2}\}$\\ [2pt]
17&$\{1|\{1|0\},\{1||1|0\}||\{1|0\},\{1||1|0\}\}$&1\\ [2pt]
18&$\{1,\{1,\{1|0\}|0,\{1,\{1|0\}|0\}\}|0,\{1,\{1|0\}|0\}\}$&$\{1|0\}$\\ [2pt]
19&$\frac{11}{16}$&$\frac{11}{16}$\\ [2pt]
20&$\{\{1||1|0||||1||1|0|||0,\{1|0\}\},\{1,\{1|\frac{1}{2}\}|\frac{1}{2}\}|0,\{1||1|0|||0,\{1|0\}\}\}$&$\{1|0\}$\\ [2pt]
\hline
\end{tabular}\caption{The heap values and their reduced canonical forms for some initial heap sizes of \textsc{GN}.}\label{tab:can_red}
\end{centering}
\end{table}

 The second column of Table \ref{tab:can_red} gives the values for heaps of sizes 1 through 20. 
  There is little evidence of any regularity in the value column. However, the group, 
G.~A.~Mesdal \cite{M} started the analysis of the value sequence of partizan subtraction games (with finite subtraction sets).\footnote{At the time of writing, \cite{FK,M} are the only two papers on the subject of partizan subtraction games.}
They showed that an approximation, the \textit{reduced canonical form}, was useful. 
In the discussion of the solutions to the problems posed in Question \ref{qu:example} we see what is lost when using this approximation. 
Table \ref{tab:can_red}, column 3, gives the reduced canonical forms, which are surprisingly simple. Table \ref{tab:fib_red} shows patterns and also hints at the relationships to the Fibonacci numbers.
 Note that in the `+0', `+1' and `+3'   columns (and others), 
the pattern alternates with even- and odd-indexed Fibonacci numbers.

\begin{table}[ht!]
\begin{center}
\begin{tabular}{c|cccccccccccccc|}
\hline    
  & +0 &    +1& +2    &+3   &+4  &+5    &+6    &+7     &+8   &+9 &+10    &+11   &+12  \\
\hline    
1    &1    \\[1pt]                                                                                            
2  &$\{1|0\}$    \\[1pt]                                                                                          
3&1/2 & 1 \\               [1pt]                                                                           
5  &$\{1|0\}$  &3/4 &$\{1|0\}$    \\ [1pt]                                                                              
8 &$\{1|1/2\}$ & 1  &$\{1|0\}$ & 5/8&  1\\[1pt]                                                                           
13 &$\{1|0\}$  &7/8 &$\{1|0\}$&$\{1|1/2\}$&1&$\{1|0\}$&11/16&$\{1|0\}$ \\ [1pt]                                                          
21&$\{1|1/2\}$  & 1&  $\{1|0\}$&$\{1|5/8\}$&1&$\{1|0\}$&13/16&$\{1|0\}$&$\{1|1/2\}$&
1&$\{1|0\}$ &21/32& 1\\                                       [1pt]
34 &$\{1|0\}$ & 15/16&$\{1|0\}$&$\{1|1/2\}$ & 1 &$\{1|0\}$ & 23/32&$\{1|0\}$&$\{1|1/2\}$ & 1&$\{1|0\}$&$\{1|5/8\}$&1\\[1pt]
55&$\{1|1/2\}$  & 1&  $\{1|0\}$&$\{1|5/8\}$&1&$\{1|0\}$&$25/32$&$\{1|0\}$&$\{1|1/2\}$&
1&$\{1|0\}$ &$\{1|21/32\}$& 1\\                                       [1pt]

\hline    
\end{tabular}
\caption{The reduced canonical forms for some initial heap sizes of \textsc{GN}.}\label{tab:fib_red}
\end{center}
\end{table}

The patterns suggested by Table \ref{tab:fib_red} relate to the Fibonacci sequence through Sturmian words.
Using these properties, we give a partition of positive integers into four subsets.
In Section \ref{sec:maintheorem}, 
 Theorem \ref{thm:main}, our main theorem, 
states that there are four `patterns' in terms of reduced canonical forms tied to the different parts.
 Its proof involves an interplay of game theory and number theory. 
Since the reduced canonical form is not widely known, we give the required game theory background 
 in Section \ref{sec:game theory}.  The number theory background is more varied, and better known, 
 thus, so as to not disrupt
 the flow, we first give the proof of Theorem \ref{thm:main} 
 in Section \ref{sec:main} and follow up with the number theory results in Section \ref{sec:fibW}. 
Section \ref{sec:disc} gives the solutions to the positions of Question \ref{qu:example} and  
discusses further directions.

\section{Main Theorem}\label{sec:maintheorem}

We are interested in the following subsets of the nonnegative integers (and we include some useful notation and number theory identities that will be explained in the sequel):
\begin{itemize}
\item  $B=\{B(n):n \geqslant 1\}=\{ \lfloor n\phi \rfloor + n : n \geqslant 1\}$,
\item  $\widehat{AB_0}=\{AB(n)+1: n \geqslant 0\}=\{2\lfloor n\phi\rfloor + n + 1 : n \geqslant 0\}$,
\item  $AB_0 = \{AB(n): n \geqslant 0\} = \{2\lfloor n\phi\rfloor + n : n \geqslant 0\}$, and $AB=\{AB(n): n \geqslant 1\}$, 
\item  $\widehat{B^2} =\{B^2(n)+1: n\geqslant 1\}=\{3\lfloor n\phi\rfloor +2n+1: n\geqslant 1\}$.
\end{itemize}
Note that, using the usual notation in this field, $B^2(n)$ means $B(B(n))$ and $AB(n) = A(B(n))$.
We define the \emph{Fibonacci sequence} by ($F_{-1}=1, F_0 = 0$) $F_1 = 1$, $F_2 = 1$, and $F_n = F_{n-1} + F_{n-2}$, $n\geqslant 3$. 

 Further, for all integers $i\ge 0$ and $n\ge 0$,  let $$G_i(n) = \lfloor i\phi\rfloor F_{2n+2} + iF_{2n+1} + F_{2n+3} -2=B^n(i)+F_{2n+3}-2.$$
 Note that $\{F_{2n+3}-2\mid n\ge 0\}= \{G_0(n)\mid n\ge 0\}$. For all $n>0$, let $G(n) = \{G_i(n)\mid i> 0\}$, and let $$\n=\widehat{B^2}\cup \{F_{2n+3}-2\mid n\ge 0\}.$$ 

\begin{lemma}[Partitioning Lemma]\label{lem:part} 
The sets $B$, $AB_0$, $\widehat{AB_0}$, and $\widehat{B^2}$ partition the nonnegative integers.
The sets $G(n)$, for $n>0$, partition the set $AB$.
\end{lemma}
\begin{table}[ht]
\begin{center}
\begin{tabular}{|l||c|c|c|c|c|c|c|c|c|c|c|c|c|c|c|}\hline
      & 0& 1& 2& 3& 4& 5& 6& 7& 8& 9&10&11&12&13&14\\ \hline\hline
$B$    & $\cdot$ & 2& 5& 7&10&13&15&18&20&23&26&28&31&34&36\\ \hline
$AB_0$    & 0& 3& 8& 11& 16& 21& 24&29&32&37&42&45&50&55&58\\ \hline
$\widehat{AB_0}$& 1& 4& 9&12&17&22&25&30&33&38&43&46&51&56&59\\ \hline
$\widehat{B^2}$& $\cdot$& 6& 14&19&27&35&40&48&53&61&69&74&82&90&95\\ \hline \hline
$G(1)$    &3& 8& 16& 21& 29& 37& 42& 50& 55& 63& 71& 76&84 & 92& 97\\ \hline
$G(2)$    &11& 24& 45& 58& 79& 100& 113& 134& 147& 168& 189& 202& 223& 244& 257\\ \hline
$G(3)$    &32& 66& 121& 155& 210& 265& 299& 354& 388& 443& 498& 532& 587& 642& 676\\ 
\end{tabular}\caption{Partitioning of heap sizes. The sets $G(n)$, for $n>1$ partition the set $AB$. The successive first differences in $G(n)$ are described by the word $W$ in Table~\ref{tab:sequences2}, with $a=F_{n+5}, b=F_{n+4}$. 
}\label{tab:partitioning}
\end{center}
\end{table}

We prove this result in two lemmas in Section \ref{sec:partitioning}; but, why this partitioning? We summarize our findings in 
the next result (only the precise behavior of (iii) is described later). Some of the values are easily described via a closed formula.
For all $n\in \Nz$, let $$\qm(n) = \frac{2}{3}\cdot\frac{4^n - 1}{4^{n}},\ \ \qp(n) = \frac{2}{3}\cdot\frac{4^{n} + \frac{1}{2}}{4^{n}}.$$ 
The sequences $\qm$ and $\qp$ approach $\frac{2}{3}$, monotonously increasing and decreasing respectively, and by inspection $\qm(0)=0, \qp(0)=1$.

\begin{theorem}[Main Theorem]\label{thm:main}
For the game \WPS:
\begin{enumerate}[(i)]
\item Let $h\in B$. Then $\rcf{h}=\cgbgame{1}{0}$.
\item Let $h\in \widehat{AB_0}$. Then $\rcf{h}=1$.
\item Let $h \in \n \setminus \{0\}$. Then  $\heap{h}\in [1/2, 1)$. Furthermore, for all $n\ge 0$, $\heap{F_{2n+3}-2} = \qm(n)$ and $\heap{F_{2n+4}-2} = \qp(n)$.
\item Let $n\in \N$ and let $h \in G(n)\subset AB_0$. Then $\rcf{h} = \cgbgame{1}{\qm(n)}$. 
\end{enumerate}
\end{theorem}
We give efficient  algorithms to find the heaps whose canonical forms are numbers in (iii) and the heaps whose reduced canonical forms are switches in~(iv). 

To illustrate the work to come we show how to compute the sequences $s(n)$ and $q(n)$, given some other assumptions.

\begin{example}
(This argument will be repeated and generalized in binary notation later.)
Observe that $\qm(n) = \frac{2}{3}\cdot\frac{4^n-1}{4^n} = (\frac{4^{n}-1}{3})/(2\cdot 4^{n-1})$. Note
that it has the odd numerator, $\frac{4^{n}-1}{3}$, and denominator $2\cdot4^{n-1}$ and so, as a fraction, 
 is in lowest terms.
Also, $\qp(n) = \frac{2}{3}\cdot\frac{4^n+\frac{1}{2}}{4^n} = \frac{2\cdot4^n+1}{3}/4^n$ and so is in lowest terms.

 Let $h=F_{2n+3}-2$  for some $n$. 
 Left can remove $F_{2n+2}$ and  leave
a heap of size $F_{2n+1}-2$. By induction,
 $\heap{F_{2n+1}-2}=\qm(n-1)$.
Right can remove $F_{2n+1}$ and leave a heap of size $F_{2n+2}-2$. 
By induction, $\heap{F_{2n+2}} = \qp(n-1)$. Since
 $\qp(n-1)-\qm(n-1) = \frac{1}{ 2^{2n-2}}$ and since 
 $\qm(n-1) =(\frac{4^{n}-1}{3})/(2^{2n-1})$
  then,
 by the simplicity rule,
 $\{\qm(n-1) \mid \qp(n-1) \}$ is of the form $\frac{a}{2^{2n}}$ for some $a$. Specifically,
 since $\qm(n-1) =  (\frac{4^{n-1}-1}{3})/(2\cdot 4^{n-2}) = (8\cdot\frac{4^{n-1}-1}{3})/4^{n}$
 then $a = 8\cdot\frac{4^{n-1}-1}{3}+1 = \frac{2}{3}(4^n-1)$. Hence $\{\qm(n-1) \mid \qp(n-1) \} = \qm(n)$.

Let $h=F_{2n+2}-2$  for some $n$. 
 Left can remove $F_{2n}$ and  leave
a heap of size $F_{2n+1}-2$. By induction,
 $\heap{F_{2n+1}-2}=\qm(n-1)$.
Right can remove $F_{2n+1}$ and leave a heap of size $F_{2n}-2$ and by induction $\heap{F_{2n}-2}=\qp(n-2)$. 
Since $\qp(n-2)-\qm(n-1) = 1/2^{2n-3}$ and  
$\qm(n-1) =(\frac{4^{n}-1}{3})/(2^{2n-1})$
  then,
 by the simplicity rule,
$\{\qm(n-1)  \mid \qp(n-2)\} =(2\frac{4^{n}-1}{3})+1)/(2\cdot2^{2n-1})=\qp(n-1)$. 
 
  Hence it suffices to show that any other options are dominated 
 or reverse out; that is, this computation gives the correct values, given this rather strong (but correct) assumption. However it will be easier to prove the properties simultaneously, for all the numbers,  which is the topic of Section~\ref{sec:numbers}.
\end{example}

\section{Combinatorial Games Background}\label{sec:game theory}
We give a very brief overview of Normal play combinatorial games, for more background  see \cite{ANW,BCG,ANS}. We then present the required 
background of reduced canonical form \cite{ANS}.

The type of \textit{combinatorial game} of interest is played by two players who move alternately; 
the game finishes after a finite sequence of moves regardless of the order of play;
there is perfect information; and there are no chance devices. The two players are called \textit{Left} and \textit{Right}.
The \textit{Fundamental Theorem of Combinatorial Games} (see \cite{ANW}[Theorem 2.1]) gives that there are four outcome classes: 
 $\mathcal{L}$---Left can force a win regardless of moving first or second; 
  $\mathcal{R}$---Right can force a win regardless of moving first or second; 
   $\mathcal{N}$---the next player to play can force a win regardless of whether it is Left or Right; 
   and  $\mathcal{P}$---the next player to play cannot force a win regardless of whether it is Left or Right.

Let $G$ be a position (in a game). The \textit{Left options} of $G$ are those positions that Left can move to (in one move) and 
the \textit{Right options} are defined analogously. Let $G^\mathcal{L}$ and $G^\mathcal{R}$ be the sets of Left and right 
options of $G$. Now $G$ can be identified with the sets of options, written $\{ G^\mathcal{L}\mid  G^\mathcal{R}\}$.

The \textit{disjunctive sum} of two positions, $G$ and $H$, written $G+H$ is the position in which a player plays in either $G$ or $H$ but not both. With the disjunctive sum as the binary operation, positions form an ordered abelian group where the partial order is: 
  $G\geqslant H$ iff $G-H$ is an $\mathcal{L}$-position; $G=H$ if and only if $G-H$ is a $\mathcal{P}$-position; 
  and $G$ is incomparable to $H$ if $G-H$ is an $\mathcal{N}$-position. With these definitions,
   Left prefers the positions that are higher in the order and Right the lower ones.
 
 Knowing the outcome class of both $G$ and $H$ is not sufficient to determine 
 the outcome class of $G+H$. A refinement is required. 

 The \textit{canonical form} of a position is obtained by eliminating \textit{dominated} options and bypassing reversible options. The canonical form can be interpreted as the number of moves advantage. For example, $\{ \mid \} = 0$ since neither player has a move,
 and $\{ 0 \mid \} = 1$ since this is the position in which Left has one move and Right none. Similarly,
 $\{ \mid 0\} = -1$, a move advantage to Right. 
 
 We recall some definitions.
 
  \begin{definition}
  The position $G$ is a number if, for all $G^L$, $G^L - G < 0$ and, for all $G^R$, $G^R - G > 0$.
  \end{definition}
  
  That is, a position is a number if every option is to a position that is worse than the original for that player. The 
  Simplicity Theorem (\cite{ANS}, Theorem~3.10, page~72) restated just for numbers will be useful for us.
  
  \begin{theorem} Let $G=\{G^L\mid G^R\}$ with $G^L$ and $G^R$ numbers. If $G^L<G^R$ then $G$ is the simplest
  number strictly between its options. Moreover, if $0\leqslant G^L<G^R\leqslant 1$ then $G=a/2^b$, where $a$ is odd
  and $b$ is the smallest integer such that $G^L<a/2^b <G^R$.  
  \end{theorem}
  
  The \textit{stops} are the best numbers that a player can obtain under alternating play.
  
 \begin{definition} The \emph{Left stop} and \emph{Right stop} of $G$, written $L(G)$ and $R(G)$ respectively, are given by
 \[L(G) =\begin{cases} G,\mbox{ if $G$ is a number;}\\
                                    \max\{R(G^L)\}, \mbox{ otherwise;}
 \end{cases}\\
R(G) =\begin{cases} G,\mbox{ if $G$ is a number;}\\
                                    \min\{L(G^R)\}, \mbox{ otherwise.}
 \end{cases}
\]
A position $G$ is \emph{hot} if $L(G)>R(G)$, and \emph{infinitesimal} if $L(G)=R(G)=0$.  \end{definition}

\subsection{Reduced Canonical Form}
\label{subsec:rcf}

\begin{definition}   Two positions, $G$ and $H$ are \emph{infinitesimally close} if  $-x < G-H < x$ for all  positive numbers $x$. We write this as $G=_I H$.
\end{definition}
  The intuition of \textit{reduced canonical form} is: if $K$ is the reduced canonical form of $G$, its game tree has the least depth of all games infinitesimally close to $G$. Another way of saying this is that $K$ is the simplest game such that $K = G + \epsilon$, for any infinitesimal $\epsilon$.  
  The concept is introduced in \cite{Cal}, however, the proof was flawed and a corrected version appears in \cite{GSi}. 
 The reduced canonical form is a relatively new tool, but shows its importance in \cite{NO,MNS,M}.
 
 The following development is taken from \cite{ANS} and the reader is invited to consult the book for the proofs of the following results.
 
 \begin{definition} We write $G\geqslant_I H$ if $G-H\geqslant -x$ for all positive numbers $x$.
 \end{definition}

Note that the relationship $\geqslant_I$ is transitive. 

\begin{definition}[Inf-reduction] Let $G$ be a position.
\begin{enumerate}
\item A Left option $G^{L_1}$ is Inf-dominated by $G^{L_2}$, if $G^{L_2}\geqslant_I G^{L_1}$.
\item A Right option $G^{R_1}$ is Inf-dominated by $G^{R_2}$, if $G^{L_1}\geqslant_I G^{L_2}$.
\item  A Left option $G^{L_1}$ is Inf-reversible (through $G^{L_1R_1}$), if  $G\geqslant_I G^{L_1R_1}$ for some
Right option $G^{L_1R_1}$.
\item  A Right option $G^{R_1}$ is Inf-reversible (through $G^{R_1L_1}$), if  $G^{R_1L_1}\geqslant_I G$ for some
Left option $G^{R_1L_1}$.
\end{enumerate}
\end{definition}

\begin{theorem}[\cite{ANS}]\label{thm:Inf0}
Assume that $G$ is not equal to a number and suppose that $G'$ is obtained by removing some Inf-dominated option (either Left or Right). Then $G=_I G'$.
\end{theorem}

\begin{lemma}[\cite{ANS}]\label{lem:Inf1}
Let $G$ be a position and suppose that $G^{L_1}$ is Inf-reversible through $G^{L_1R_1}$.
Let $$G' = \{G^{L_1R_1L},G^{L'}\mid G^\mathcal{R}\}$$
where $G^{L_1R_1L}$ ranges over all Left options of $G^{L_1R_1}$ and $G^{L'}$ ranges over all 
Left options of $G$ except $G^{L_1}$. If $G'$ is not a number then $G=_IG'$.
\end{lemma}

\begin{theorem}[\cite{ANS}]\label{thm:Inf2}
Assume that $G$ is hot and suppose that some Left option $G^{L_1}$ is Inf-reversible through $G^{L_1R_1}$.
Let $G'$ be as in Lemma \ref{lem:Inf1}. Then $G=_IG'$.
\end{theorem}

A position $H$ is a \emph{sub-position} of a position $G$ if there is a (not necessarily alternating and possibly empty) sequence of consecutive moves from $G$ to $H$.

\begin{definition} A position $K$ is in \emph{reduced canonical form}, if, for every sub-position $H$ of $K$, either $H$ is in canonical form and is a number; or $H$ is hot and $H$ does not contain any Inf-dominated or Inf-reversible options.
\end{definition}

The next result is a combination of results in \cite{ANS}.

\begin{theorem}
For every position $G$ there exists a unique reduced canonical form position $K$ such that $G=_I K$.
\end{theorem}

\begin{definition}
  The \emph{reduced canonical form} of a position $G$, denoted $\rcf{G}$, is the unique reduced canonical form position, $K$, such that $\rcf{G} = K$.
\end{definition}

The next result will be useful in this paper since showing $G\geq_I K$ is equivalent to showing $G-K\geq_I0$.
\begin{lemma} The following are equivalent.
\begin{itemize}
\item $G\geq_I 0$;
\item $R(G)\geq 0$;
\item $G\geq \epsilon$ for some infinitesimal $\epsilon$.
\end{itemize}
\end{lemma}

\begin{corollary}\label{cor:numberish}
If $L(G)=R(G)=x$ then $G=x+\epsilon$ for some infinitesimal $\epsilon$.
\end{corollary}
\begin{proof} Since $L(G)=R(G)=x$ then $x$ is a number. By the Number Avoidance Theorem, $L(G-x)=R(G-x) = x-x=0$
and thus $G-x$ is an infinitesimal.
\end{proof}

\begin{example}
The position $1$ is incomparable with $\cgbgame{1}{0}$, 
but $1 + \cgstar + \cgbgame{0}{-1} > 0$ so we conclude that $1 \geq_{I} \cgbgame{1}{0}$.
\end{example}

Let $a$ and $b$ be numbers. Then 
$\cgbgame{a}{b}$ is a \emph{switch} if $a \geqslant b$ (and otherwise a number). 
Let $G=\cgbgame{a}{b}$ be a switch and let $c$ be a number. We have the following comparisons:
if $c > a$, then $c > G$; if $a \geqslant c \geqslant b$, then $G$ is incomparable with $ c$; if $b >c$, then $G > c$.

\begin{lemma}\label{lem:bigish0}
 If $x$ is a number and $x\geqslant 0$, then $\cgbgame{x}{0} \geqslant_I 0$.
\end{lemma}

\begin{proof}
If $x=0$ then $\{0\mid 0\} =_I0$. Suppose $x>0$. Since $\cgbgame{x}{0} + \cgup >0$ 
we know that $\cgbgame{x}{0} \geqslant_I 0$. 
\end{proof}
Note that since $\cgbgame{x}{0}$ is in reduced canonical form, and so is $0$, but they are not equal, so $\cgbgame{x}{0} \not =_I 0$ 
and we can also conclude that $\cgbgame{x}{0} >_I 0$, that is $\rcf{\{x\mid 0\}}>0$.

\section{\WPS}\label{sec:main}

We provide a general bound for any \WPS position. 

Note that in general, a game position and its value are not usually distinguished but since the heaps sizes and each value is a number 
or a switch of numbers, it is convenient to retain the $\heap{h}$ notation. We abuse the convention and sometimes say
`the move to $\heap{h}$' as shorthand for `a move to heap of size $h$ which has value $\heap{h}$'. This will be clear in context.

\begin{lemma}\label{lem:rcfbound}
 For all $h >0$, $\cgbgame{\frac{1}{2}}{0} \leqslant \rcf{h} \leqslant 1$.
\end{lemma}

\begin{proof}
The result holds for heap size 1. For any $h>1$, 
$\heap{h} = \cgbgame{\heap{h-A(i)}}{\heap{h-B(i)}}$ for $i>0$ and $A(i),B(i)\leqslant h$. By induction, 
$$\heap{h} \leqslant_I \cgbgame{1+\epsilon_i}{1+\epsilon_j} =1+\{\epsilon_i\mid \epsilon_j\}=_I 1$$
where $\epsilon_i$ and $\epsilon_j$  are infinitesimals (see Corollary \ref{cor:numberish}).

Note that $\heap{0} = 0$, and that the result also holds for heap sizes up to and including 3. The Left and Right stops of 
$\heap{4} = \{1 \mid\mid 1\mid 0\}$ are 1; thus $\heap{4}=1+\epsilon =_I 1$; thus the result is true up to $h=4$.

 From any larger heap, Left has an option to at least one of $\heap{3} = \frac{1}{2}$
  or $\heap{4} =_I 1$. By induction, $\heap{h} \geqslant \cgbgame{\frac{1}{2}}{\epsilon} =_I \cgbgame{\frac{1}{2}}{0}$. 
  That is, $\rcf{h}\geqslant\{\frac{1}{2}\mid 0\}$.\end{proof}
  
  In the next subsections  we prove the items of the main theorem.

\subsection{Theorem~\ref{thm:main} (i)}

In this subsection, we consider $h \in B$, the only case where Right has a winning move. 

\begin{theorem}\label{lem:1slash0}
If $h \in B$, then $\rcf{h} = \cgbgame{1}{0}$.
\end{theorem}
\begin{proof}
Right has a move to $0$, which dominates (or Inf-dominates) any other option. By Lemma \ref{lem:rcfbound} and 
because $0 <_I \cgbgame{\frac{1}{2}}{0}$ by Lemma \ref{lem:bigish0}. Since $h \in B$, then $h-1 \in A$, 
so Left has a move to $1$, and $\heap{1} = 1$ which Inf-dominates any other option, by Lemma \ref{lem:rcfbound}.
\end{proof}

\begin{corollary}\label{cor:not to B}
A move to a heap-size in $B$ is Inf-reversible to 0 for Left and to 1 for Right.
\end{corollary}
\begin{proof}
From Theorem \ref{lem:1slash0}, we know that for $h \in B$, $\heap{h}$ has reduced canonical form 
$\cgbgame{1}{0}$. From Lemma \ref{lem:rcfbound}, we can see that the move to $\heap{h}$ is Inf-reversible.
\end{proof}
The consequence of this result is that only rarely does a player want to move to a heap size in $B$ (see also Section~\ref{sec:disc}).

\subsection{Theorem~\ref{thm:main} (ii)}

Recall $$\widehat{AB_0}=\{AB(n)+1: n \geqslant 0\}=\{2\lfloor n\phi\rfloor + n + 1 : n \geqslant 0\}$$
and we prove that these heap sizes all have reduced canonical form 1. 

\begin{theorem}\label{thm:rcf1}
 If $h \in \widehat{AB_0}$, then $\rcf{h} = 1$.
\end{theorem}

\begin{proof}
Let  $h = AB(n) + 1$ for some $n \geqslant 0$. The result is true for $n = 0$ so we assume $n > 0$. 

Left, going first, can move to $\heap{1} = 1$ by removing $AB(n)$. 
If Right can move to a position where Left cannot move to $\heap{1}$ then Right can move to $\heap{1}$ from that position; thus there exist positive $i$ and $j$ such that $AB(n) + 1 - B(i) - B(j) = 1$, 
which contradicts Lemma \ref{lem:bbnotabn}.
Thus $L(\heap{h})=R(\heap{h})=1$ and, by Corollary \ref{cor:numberish}, $\rcf{h} = 1$.
\end{proof}

\subsection{Theorem~\ref{thm:main} (iii)}\label{sec:numbers}

In this section we are concerned with the canonical forms of the heaps in $\n$ and we show that they are numbers.
The \textit{even} (or ternary) \textit{Fibonacci representation} of a nonnegative integer, defined in Section \ref{sec:fibW}, is using only the even indexed Fibonacci numbers and it is unique. We will show what the value is when a player removes the largest Fibonacci number available and then the 
`fundamental Claim for the numbers' (backed up by some number theory on the even Fibonacci representation)
shows that all other possible moves give a value that is dominated or reverses out. Thus the value of the original position
is $\{a\mid b\}$ for some numbers $a<b$ and the actual value is given by the standard CGT Simplicity Theorem.
In this case, they will be obtained by taking the mean of numbers in smaller heaps. Suppose that $x<y$ 
then, the simplest number between the options in $\cgbgame{x}{y}$ is the value of this position and the depth of the canonical form tree is the length of the binary representation of the number. 

In addition to the Zeckendorf representation and the even (Fibonacci) representation, in this section we also use the standard binary representation of rational numbers. In binary notation, for $n>0$, $$\qm (n) = 0.(10)^{n-1}1,\  \ \qp (n) = 0.(10)^{n-1}11.$$
In fact, if we wish to extend the notation to all $n\ge 0$, then we define $\qm(n)=(01)^n0$ and $\qp (n) = (01)^n1$, and insert the `$.$' as appropriate for $n>0$.
The consistency with the definition of the sequences is easily justified via partial sums of these geometric series.

\begin{table}[tbp]
  \centering
  \begin{tabular}{|c|clc|ccc|} \hline
heap&value&bin.repr.&obtained from&optimal moves&options&option's values\\\hline
0&0&0&&&&\\
1&1&1&&&&\\\hline
\bf{3}&$\frac{1}{2}$&0.1&$\cgbgame{0}{1}$&3,2&0,1&$\cgbgame{0}{1}$\\[3pt]\hline
\it{6}&$\frac{3}{4}$&0.11&$\cgbgame{0.1}{1}$&3,5&3,1&$\cgbgame{0.1}{1}$\\[3pt]\hline
\bf{11}&$\frac{5}{8}$&0.101&$\cgbgame{0.1}{0.11}$&8,5&3,6&$\cgbgame{0.1}{0.11}$\\[3pt]\hline 
14&$\frac{7}{8}$&0.111&$\cgbgame{0.11}{1}$&8,13&6,1&$\cgbgame{0.11}{1}$\\[3pt]
\it{19}&$\frac{11}{16}$&0.1011&$\cgbgame{0.101}{0.11}$&8,13&11,6&$\cgbgame{0.101}{0.11}$\\ [3pt]\hline
27&$\frac{13}{16}$&0.1101&$\cgbgame{0.11}{0.111}$&21,13&6,14&$\cgbgame{0.11}{0.111}$\\[3pt]
\bf{32}&$\frac{21}{32}$&0.10101&$\cgbgame{0.101}{0.1011}$&21,13&11,19&$\cgbgame{0.101}{0.1011}$\\[3pt]\hline
35&$\frac{15}{16}$&0.1111&$\cgbgame{0.111}{1}$&21,34&14,1&$\cgbgame{0.111}{1}$\\[3pt]
40&$\frac{23}{32}$&0.10111&$\cgbgame{0.1011}{0.11}$&21,34&19,6&$\cgbgame{0.1011}{0.11}$\\[3pt]
48&$\frac{27}{32}$&0.11011&$\cgbgame{0.1101}{0.111}$&21,34&27,14&$\cgbgame{0.1101}{0.111}$\\[3pt]
\it{53}&$\frac{43}{64}$&0.101011&$\cgbgame{0.10101}{0.1011}$&21,34&32,19&$\cgbgame{0.10101}{0.1011}$\\[3pt]\hline
61&$\frac{25}{32}$&0.11001&$\cgbgame{0.11}{0.1101}$&55,34&6,27&$\cgbgame{0.11}{0.1101}$\\[3pt]
69&$\frac{29}{32}$&0.11101&$\cgbgame{0.111}{0.1111}$&55,34&14,35&$\cgbgame{0.111}{0.1111}$\\[3pt]
74&$\frac{45}{64}$&0.101101&$\cgbgame{0.1011}{0.10111}$&55,34&19,40&$\cgbgame{0.1011}{0.10111}$\\[3pt]
82&$\frac{53}{64}$&0.110101&$\cgbgame{0.1101}{0.11011}$&55,34&27,48&$\cgbgame{0.1101}{0.11011}$\\[3pt]
\bf{87}&$\frac{85}{128}$&0.1010101&$\cgbgame{0.10101}{0.101011}$&55,34&32,53&$\cgbgame{0.10101}{0.101011}$\\[3pt]\hline
\end{tabular}  
  \caption{We display heaps that are numbers and their binary representations; the optimal moves from the respective heaps, the corresponding options and their values. The {\bf bold} and {\it italicized} heaps correspond to $\heap{F_{2n+3}-2} = \qm(n)$ and $\heap {F_{2n+4}-2}=\qp(n)$, respectively.}
  \label{table:numbers}
\end{table}
From the unique finite binary representation $d=d_0.d_1\ldots d_{k}$ of a dyadic number $1/2 \leqslant d \leqslant 1$ we assign a positive integer 
\begin{align}\label{xirepr}
\xi=\xi(d)=\sum d_iF_{e(i)}
\end{align} 
where $e(i)$ is a function defined recursively from the digits of $d$ to even positive integers: $e(0)=2, e(1)=4$,
\begin{equation}\label{e}
e(i)=\left\{
\begin{array}{ll}
e(i-1)  & \text{ if } d_{i-2}d_{i-1} = 01,\\ 
e(i-1)+2,  & \text{ otherwise}. 
\end{array}\right.
\end{equation}

For example $d = 0.110011$ will produce $\xi(d)=116$, namely $e(0)=2$, $e(1)=e(2)=4, e(3)=e(2)+2=6, e(4)=e(3)+2=8, e(5)=e(4)+2=10,e(6)=e(5)=10$, so that $\xi=F_4+F_4+F_{10}+F_{10}=3+3+55+55$.

\begin{lemma}\label{lem:even}
The Fibonacci representation of $\xi$ in (\ref{xirepr}) is the even representation.
\end{lemma}
\begin{proof}
By definition of $e$, it is clear that only even-indexed Fibonacci numbers are included in the Fibonacci representation of $\xi$. 
By definition, the only way that $e(i)=e(i-1)$ is if $d_{i-2}d_{i-1} = 01$, which in its turn implies that $e(i)<e(i+1)$. 
Hence the representation is (at most) ternary. 
If both $2F_{2i}$ and $2F_{2j}$, $i>j$,  occur then they correspond to different strings `011' in $d$, at $d_{i-2}d_{i-1}d_i$ and $d_{j-2}d_{j-1}d_j$ respectively; 
thus there is a leftmost intervening pattern $d_{k-2}d_{k-1}d_k=110$ in $d$. By (\ref{e}), this implies that $F_{2k}$, $i>k>j$, does not occur in the representation. \end{proof}

We wish to convert any heap size $h$ in our list $\n$ to a unique number in binary via $\xi^{-1}$ and then show that this number is 
$\heap{h}$. For example, $h=116$ is easily seen to produce the above binary number; the two 2s correspond to two factors
 $011$ and there should be two even Fibonacci numbers missing so we must include another $0$ between these factors. 

\begin{lemma}\label{heap_number}
Consider a heap size $h\in \n \setminus \{0\}$. Then there is a unique number $1/2\le \delta \le 1$ in binary such that $\xi (\delta)=h$.
\end{lemma}
\begin{proof}
Let $h=\sum c_iF_{2i}$, where $c_i\in \{0,1,2\}$, be the even Fibonacci representation. We need to define $\delta$ recursively, given the function $e$, such that $h=\sum \delta_jF_{e_j}$ is the even representation (by Lemma~\ref{lem:even}). In order to know $e$, we need exactly two preceding bits of $\delta$. Hence let us begin by finding $\delta_0\delta_1$.

By previous results we know that $F_4$ is the least representative, which means that it is given that $c_0=0$.
That each integer in $\n$ has at least one $F_4$ is clear if it is of the form $F_n-2, n\geqslant 5$.  It follows from Lemma~\ref{lem:repr};
 if it is of the form $B^2(n)+1$ with $n>0$. Thus $c_1\in\{1,2\}$.
 Hence the prefix of $\delta$ is $0.1$, which coincides with the given bounds $1/2 \le \delta \le 1$. 

We divide the rest of the proof into cases. Suppose that $ c_{i-1}$, in the even representation of $h$, is obtained from 
$\delta_{j-2}\delta_{j-1}$. Suppose $\delta_{j-2}\delta_{j-1}=01$, with $j\ge2$. It follows that
$c_{i-1}\ne 0$ and $e(j)=e(j-1)$. If $c_{i-1}=1$, then we let $\delta_j=0$, and otherwise $\delta_j=1$. 

It remains to check the cases $\delta_{j-2}\delta_{j-1}\ne 01$. They are similar in the sense that 
\begin{align}\label{ej2}
e_j = e_{j-1}+2. 
\end{align}
Hence, we only influence the number of $F_{2i}$s (the number of  $F_{2i-2}$s is assumed to have been translated correctly by induction). 
\begin{itemize}
\item If $c_i=0$, then we must let $\delta_j=0$. This suffices, since, by (\ref{ej2}), $F_{2i} = F_{e(j)}$ will correctly be omitted, and, by (\ref{e}), which gives $e(j+1)=e(j)+2$. 
\item If $c_i=1$, then we must let $\delta_j=1$. By (\ref{ej2}), $F_{2i} = F_{e(j)}$  but we do \emph{not} want another copy. 
Hence we must look into the three cases $\delta_{j-2}\delta_{j-1}\delta_j = 001, 101 \text{ or } 111$. The last case will give 
$e(j+1)=e(j)+2$, so this will be correct by default. For the first two cases, we need to put $\delta_{j+1}= 0$, because 
otherwise $e(j)=e(j+1)$ would give two copies of $F_{e(j)}$. Then $\delta_j\delta_{j+1}= 10$ will give $e(j) < e(j+1)$.
Thus each instance gives a unique correct update of $\delta$. 
\item If $c_i=2$. This case is similar to the previous one. We must let $\delta_j=1$. By (\ref{ej2}), $F_{2i} = F_{e(j)}$ will correctly be included, and we do want  another copy. Hence we must look into the three cases $\delta_{j-2}\delta_{j-1}\delta_j = 001, 101 \text{ or } 111$. For the first two cases, we need to put $\delta_{j+1}= 1$, because then (and only then) $e(j)=e(j+1)$ gives two copies of $F_{e(j)}$.  The last case will give $e(j+1)=e(j)+2$, and so there is no way to produce another copy of $F_{e(j)}$. However, 
by the even representation and by definition of the function $e$ this (the only remaining) case cannot happen. Specifically, because
 there has to be a rightmost 0 to the left of a leftmost 1 in this factor of consecutive 1s. By the definition of $e$, 
 this is an instance of a 2 of an even indexed Fibonacci number, say $2F_{2k}$, with $k<j$. Since there are only 1s in the factor, each even indexed Fibonacci number between $F_{2k}$ and $F_{2j}$ is included. Then by the even representation, $F_{2j}$ can be included at most once.
\end{itemize}
\end{proof}

By this result, for each heap-size $h\in \n$, $\xi^{-1}$ finds, in polynomial time, the dyadic rational $\delta(h)$. 
It remains to prove that, in fact, $\delta(h)=\heap{h}$. 

Recall that Left removes even indexed and Right removes odd indexed Fibonacci numbers.
 Let us denote by $F_{{\rm Rmax}}$ the largest Fibonacci number that Right can remove (odd indexed), and by $F_{{\rm Lmax}}$ the largest Fibonacci number that Left can remove (even indexed), given a position $h\in \n$.

\begin{proposition}\label{lem:maxFib}
The $\xi$-algorithm is consistent with the evaluation of the values for a heap $h\in \n$ in the following sense. \begin{itemize}
\item[{\rm (i)}] If $F_{{\rm Lmax}} > F_{{\rm Rmax}}$ and  $\heap h= x01$ (a binary fraction) then $\xi (x01)=h$, and $\xi (x)=h-F_{{\rm Lmax}} < h-F_{{\rm Rmax}}=\xi(x1)$.
\item[{\rm (ii)}] If $F_{{\rm Rmax}} > F_{{\rm Lmax}}$  and  $\heap h= x01^r$, with $r\in \N$, then $\xi (x01^r) = h$, and $\xi (x01^{r-1})=h-F_{{\rm Lmax}} > h-F_{{\rm Rmax}}=\xi(x1)$.
\end{itemize}
\end{proposition}
\begin{proof}
Recall that the Simplicity Theorem for Games gives that if $G$ is a position and $G^L<G^R$ and both are numbers then $G$ is a number
and specifically, if $x$ is a binary fraction between 0 and 1 then $x01^r=\cgbgame{x01^{r-1}}{x1}$ for $r\geq 1$.

We begin with $x$, the empty (binary) word, and interpret the resulting dyadic rational as 0.1 = 1/2, and $\xi(0.1)=F_4=3$, by (i). The number is the simplest with denominator 2 and it is the mean of 0 and 1. If $r=2$ and $x$ is empty, then we claim that the heap size is $\xi(0.11)=2F_4 = 6$. We apply (ii), since $3/4=0.11 = \cgbgame{0.1}{1} = \cgbgame{1/2}{1}=\cgbgame{\heap 3}{\heap 1}$, and $1<3$. The next number is obtained by letting $x = 0.1$. We get $h=F_4+F_6=11$ and, by (i), $\xi^{-1}(11)=0.101=\cgbgame{0.10}{0.11}=\cgbgame{1/2}{3/4}=5/8$. We can think of this procedure as an algorithm, by noting that, given that all values are numbers, we must have $h^L < h^R$. Claim: by recursively applying any thus legal combination of numbers, we obtain all \wps\ numbers, and by applying $\xi$, we get each heap size that is a canonical form number. (Given a heap size in $\N$ there is a much faster way to find out its value, but it is irrelevant for this proof.)

By the binary notation, it is immediate that the value computation has smallest possible denominator 
(which is $2^n$, where $n$ is the index of the rightmost ``1''), and also that in both cases the value is the mean of the options. 
By the simplicity theorem, this part is correct. 
It remains to verify that the prescribed value options are actual heap size options, and also that they are optimal. They are actual heap size options, because, in case (i), by definition of $\xi$, by going from $x01$ to $x$, the largest even indexed Fibonacci option has been removed from $\xi(x01)$. This is an option for Left, by number theory section. Similarly, $\xi(x01)-\xi(x1)=F_{2n+4}-F_{2n+2} = F_{2n+3}$, which is odd indexed, hence it is in $B$,  which is Right's subtraction set. In case (ii), Left's option is as in case (i), the largest even indexed Fibonacci number, as defined by $\xi$, is removed. For Right,  if $r>1$ then the situation is somewhat different. If $r=2$, then, for some index $n$, we compute  $2F_{2n}-F_{2n-2}=F_{2n+1}$. Otherwise the pattern $01^r$ (in the value) corresponds, via $\xi$, to $2F_{2n}+F_{2n+2}+\cdots +F_{2n+2(r-2)}$. Now $2F_{2n}+F_{2n+2}+\cdots +F_{2n+2(r-2)}-F_{2n-2}=F_{2n+2r-1}$, which is the largest odd indexed Fibonacci number smaller than $\xi(x01^r)$.

Note also that in (i) and (ii), the inequalities of the maximal Fibonacci numbers are consistent with the definition of $\xi$.

The result follows by induction by proving the optimality of the given options. By induction and previous results, we only need to consider the heaps in $\n$. It remains to prove the following claim.\\

\noindent \textbf{The Fundamental Claim of the \wps \ numbers:} Let $h\in N$. Let $h^L$ and $h^R$ denote the option, where Left and Right has subtracted the largest available Fibonacci number, respectively.
Then, for each $h>x\in \n$, $ d(x) > d(h^L)$ iff $h-x\in B$ and  $d(x) < d(h^R)$ iff $h-x\in A$.  
\\

This Claim is better proved in purely number theoretic terms in the following equivalent statement. 
Let $z_1(x)$ denote the smallest term in the Zeckendorf representation of $x$.
\begin{theorem}\label{thm:new}
Let $d, g\in [1/2,1)$, such that $\xi(d), \xi(g)\in \n$. Suppose that $s:=\xi(d) - \xi(g)>0$. Then $z_1(s)$ is \emph{odd} if and only if $g>d$.
\end{theorem}

The idea is that the position of value $g$ is a Right option from the position of value $d$ if and only if $g > d$ 
(but Right does not want to play there). The position of value $g$ is a Left option if and only if $g < d$ 
(but Left does not want to play there). This explains why the $\xi$-algorithm gives values that are {\it numbers}. 
Note however that the proof of the theorem does not require any `game reasoning'--- it is pure number theory. 
We keep this result in the game section, because it concerns the $\xi$-algorithm, which is introduced here.

\begin{proof}[Proof of Theorem \ref{thm:new}]
By Lemma~\ref{heap_number}, $g \ne d$. Let $k$ be the index of the most significant bit where the binary representations of $g$ and $d$ differ. Suppose that $g>d$. Then 
\begin{align}\label{gk>dk}
1=g_k > d_k=0. 
\end{align}
Since $d_k=0$ and $k\ne 0$ (since $g_k=1$), we know that $k>1$. By $\xi$, this means that the least index in the even representation of $\xi (g)$ is at least 4. By (\ref{gk>dk}), the least index in the even representation of $\xi (d)$ is at least 6. Let us write the difference $\xi(d)-\xi(g)$ in the even representation as 
\begin{align}\label{stuff}
\xi(d)-\xi(g) = \sum d_iF_{e(i)} - \sum g_iF_{e(i)} = \sum \eta_i F_{2i}, 
\end{align}
where $\eta_i \in \{-2,-1,0,1,2\}$.

Algorithm: Let $\gamma_+ = \max \{i\mid \eta_i\ne 0\}$. Then $\eta_{\gamma_+}>0$. Let 
\begin{align}\label{max-}
\gamma_- = \max \{i\mid \eta_i < 0\}.
\end{align}
We want to successively subtract each negative Fibonacci term, by starting with the one of largest absolute value 
$$\sum_{\eta_i>0} \eta_i F_{2i}\ -\ F_{2\gamma_-}.$$
By each transformation, we will write each large Fibonacci term in the Zeckendorf representation (and verify that, for the next step, we can ignore all but the smallest of the thus obtained larger terms). The small (positive) terms will remain in the even representation.
We will repeat this algorithm until the definition of $\gamma_-$ gives no output. Then we want to conclude that the obtained Zeckendorf representation of (\ref{stuff}) has a least odd index.

\emph{First step}: Let $2j$ be the smallest index of a positive Fibonacci term greater than $F_{2\gamma_-} $. The first step gives  $F_{2j} - F_{2\gamma_-} = F_{2i-1} +F_{2i-3} +\cdots + F_{2\gamma_-+1}$, and we observe that $F_{2\gamma_-+1} > F_{2\gamma_-} > F_{2\gamma'_-}$, where $2\gamma'_-$ is the new index obtained via (\ref{max-}), by ignoring $\eta_{\gamma_-}$. Let us study the case 
\begin{align}\label{eta2}
\eta_{\gamma'_-}=2
\end{align}
 (otherwise the first step is done). In this case, we rewrite the term $2F_{2\gamma'_-} = F_{2\gamma'_- + 1} + F_{2\gamma'_- -2}$. We make two observations:
\begin{itemize}
\item[(i)] The largest remaining Fibonacci term to subtract is instead $F_{2\gamma'_- +1}$
\item[(ii)] The term $F_{2\gamma'_- -2}$ will be the largest subtraction term in the next step, and perhaps it will be $2F_{2\gamma'_- -2}$.
 Note also that, by the even representation and (\ref{eta2}), it cannot be $3F_{2\gamma'_- -2}$. 
\end{itemize}
The subtraction in (i) will give at most $F_{2\gamma'_-}$ as the currently smallest term in what is to become the Zeckendorf representation of $s$. Then, by (ii), the algorithm is correct; the remaining two cases being as follows.  

For the \emph{second step} (and onwards), there is a possibility that the index  is odd: if so, let $2j+1$ be  the smallest index of a positive Fibonacci term greater than $F_{2\gamma'_-} $.

This leads to 
\begin{align}\label{nothingnew}
F_{2j+1} - F_{2\gamma'_-} = F_{2j} +F_{2j-2} +\cdots + F_{2\gamma'_-+2}+F_{2\gamma'_- -1}.
\end{align}
This time, we observe that $F_{2\gamma'_- -1} > F_{2\gamma''_-}$, where $2\gamma''_-$ is the new index obtained via (\ref{max-}), by ignoring $\eta_{\gamma_-}$ and $\eta_{\gamma'_-}$. Perhaps $F_{2\gamma'_- -1} = F_{2\gamma''_- +1}$. In this case $F_{2\gamma'_- -1} - F_{2\gamma''_-} = F_{2\gamma''_- -1}$, which (with respect to our algorithm) is analogous to (\ref{nothingnew}). It remains to discuss the case as in (\ref{eta2}). It will result in a subtraction of the form $F_{2j+1} - F_{2\gamma'_- +1}$. This could equal zero, or otherwise it is a sum of consecutive even indexed Fibonacci numbers, the smallest index being $2\gamma'+2$, which is at least two larger than the largest remaining index of negative Fibonacci terms.

This algorithm can terminate in two different ways. Either $\eta_- = 1$ or $\eta_- = 2$ (where the previous negative terms have been omitted). In the first case we subtract an even indexed Fibonacci number from another (odd or even indexed) Fibonacci number. This results in an odd indexed Fibonacci number (with index $\ge 5$, since the subtracted term has index at least 4). In the second case, we apply a rewrite as for (\ref{eta2}). The extremal case of $2F_4$ cannot happen, because the smallest possible index for a Fibonacci term to subtract is 6. By applying the algorithm, this could lead to at most one single $F_4$, that is: if applying (\ref{eta2}) gives $2\gamma =4$, then $\eta_- = 1$. Otherwise we rewrite $2F_{2\gamma}=F_{2\gamma -2}+F_{2\gamma -1}$, and we can apply algorithm with $\eta_-=1$ and $2\gamma \ge 4$.

Therefore we get that $z_1(s)$ is odd, which concludes this part of the proof.\\

Suppose next that $g < d$. This implies that there is a smallest index $k>1$, such that $g_k\ne d_k$, and then $g_k=0$ and $d_k=1$. We must show that $z_1(s)$ is even. Suppose first that $k$ is the only index with differing bits. If the $1=d_k$ introduces a new 01 factor, then, by $\xi$, $s<0$, so this is impossible. Hence, all other Fibonacci terms must be identical, and so $s=F_{2i}$, some $i$, since $\xi$ gives the even representation.

Suppose next, that there is more than one differing bit. If each differing bit $d_i\ne g_i$ satisfies $d_i=1$, then, by $s$ not in the form $F_2+\cdots +F_{2i}+2F_{2i+2}$, $i\ge 0$, we use Lemma \ref{evenreprZodd} to conclude that $z_1(s)$ is even.

The next case is that there is an index $i>k$ that satisfies $1=g_i>d_i = 0$. In this case, there must be at least one more differing bit of index $j > i$ such that $0 = g_j< d_j = 1$, for otherwise $s\le 0$. 
We apply the algorithm in the first part of the proof to the larger terms in absolute value.
This will result in the least odd indexed term in the Zeckendorf representation having a larger index than that of
 the even indexed Fibonacci term resulting from the $k$th bit. It is then an easy computation to see that $z_1(s)$ is even.
\end{proof}

This finishes the proof of Proposition \ref{lem:maxFib}.
\end{proof}
Thus we have a proof for the main result for the numbers of \WPS.
\begin{theorem}
For all $h\in \n$, $d(h)=\heap{h}$.
\end{theorem}
\begin{proof}
This follows by Proposition \ref{lem:maxFib}.
\end{proof}
Note that it follows by the $\xi$ algorithm (Lemma~\ref{heap_number}) that the \WPS numbers are pairwise distinct.

\subsection{Theorem~\ref{thm:main} (iv)}
From Section \ref{sec:numbers},
 each heap size of the form $G_0(n) = F_{2n+3}-2$, for $n\in \Nz$,  has value $\qm(n)$. 
 Here we prove that, for $i,n\in \N$, a heap size of the form $G_i(n) = F_{2n+2}\lfloor i\phi\rfloor + F_{2n+1}i + F_{2n+3} -2$ 
 has reduced canonical form
$\cgbgame{1}{\qm(n)}$. Note here that both $i, n > 0$.

There is a simple proof for the characterization when the reduced canonical form is $\cgbgame{1}{\frac{1}{2}}$. The rest of the result is contained in the somewhat-more-challenging-to-prove Theorem \ref{thm:Q-}.

\begin{proposition}\label{prop:1slash1/2}
If $h \in \{3\lfloor n\phi\rfloor + 2n + 3 \mid\ n \geqslant 1\}$, then $\rcf{h} = \cgbgame{1}{\frac{1}{2}}$.
\end{proposition}

\begin{proof}By equation (20), $A(n)+2B(n) = B^2(n)$ and hence $3\lfloor n\phi\rfloor + 2n \in B$. 
Since Right has the option to move to heap size 3, ($\cf{3} = \frac{1}{2}$) and cannot move to $0$, then $RS(h) = \frac{1}{2}$. 
Since no two consecutive integers are in $B$, then Left
has the option to move from $h$ to heap size 4 and consequently, $LS(h)=1>\frac{1}{2} = RS(h)$ and so $h$ is hot.
Thus Theorems \ref{thm:Inf0} and \ref{thm:Inf2} allow us to replace every option by its reduced canonical form.
Consequently, $\rcf{h} = \{1\mid 1/2\}$.
\end{proof}

\begin{theorem}\label{thm:Q-}
If $h\in G(n)$, then $\rcf{h}=\cgbgame{1}{\qm(n)}$.
\end{theorem}

\begin{proof}
We will prove that $h$ is hot by showing: (i) Left has an option to a heap size in $AB_0$ and thus $LS(h) =1$; 
(ii) Right has a move to a position with value $s(n)$ and thus $RS(h)<1$. Therefore, 
Theorems \ref{thm:Inf0} and \ref{thm:Inf2} allow us to replace every option by its reduced canonical form.

By assumption, we let $h\in \{AB(i)\}\setminus \{F_{2n+3}-2\}$ and wish to show that $\rcf{h} = \cgbgame{1}{\qm(n)}$. By Theorem~\ref{theorem:genK} we have $\lfloor i\phi\rfloor  F_{2n+2}+ iF_{2n+1} = F_{2n}A(i) + F_{2n+1}B(i) = B^{n+1}(i)\in B$, for all $n\geqslant 0$ and all $i>0$. Hence it is clear that Right can move to position $F_{2n+3} -2$. By induction, the value of this position is $\qm(n)$. Since $\qm$ is strictly increasing, and any other number is greater than $2/3$, it suffices to  show that Right cannot move to any $\qm(m)$ for $m<n$. 

Note that, for $i=1$, by Corollary \ref{corollary:fibadd}, the difference 
\begin{align*}
x=x(m,n) &:=F_{2n+3}-2 - (F_{2m+3}-2)\\ &= F_{2n+3} - F_{2m+3}\\ &=\sum_{m+2\leqslant j\leqslant n+1}F_{2j},
\end{align*}
 is in $A$, for all $0\leqslant m<n$. Buy Theorem~\ref{thm:RA}, it now follows that $x + F_{2n+3}\in A$, since $x\leqslant F_{2n+3}+1$. Hence Right cannot move from position $2F_{2n+3}-2$ to position $F_{2m+3}-2$, for any $m,n$.
 
For the general proof, note that we just showed that both 
\begin{align*}
x_0=x_0(m,n) &:= F_{2n+3}-2 - (F_{2m+3}-2)\\ &= F_{2n+3} - F_{2m+3} 
\end{align*}
and 
\begin{align*}
x_1=x_1(m,n) &:= 2F_{2n+3}-2 - (F_{2m+3}-2)\\ &= 2F_{2n+3} - F_{2m+3}\\&=x_0+F_{2n+3}, 
\end{align*}
are in $A$, for any $n,m \geqslant 0$. We will lift this result recursively to all $i$. Note that, for $i=2$, we can define $x_2=x_1+2F_{2n+2}+F_{2n+1}=x_0+F_{2n+5}$ and where $x_0< F_{2n+5}$. Hence, by Theorem~\ref{thm:RA}, it also follows that $x_2\in A$. 

The result follows from the following number theory lemma (we keep it here, because it is tightly related to the game strategies). The cases $x_0, x_1$ and $x_2$ have already been checked. We will often use that $(x_i)$ is increasing. 

\begin{lemma}\label{lem:switch}
Let $n \in \Nz$, and let $x_1=x_0+F_{2n+3}$, where $x_0=F_{2n+3}-F_{2m+3}< F_{2n+3}$, for some $m<n$. Further, for any integer $k\ge 3$, and $j = 0, \ldots , F_{k-1}-1$, let 
\begin{align}\label{x}
x_{F_k+1+j}=x_{j+1}+F_{2n+k+2}.
\end{align}
Then, for all $k, n$, $x_{j+1} < F_{2n+k+2}$ and , for all $i\in \Nz$, $x_i=G_i(n)-F_{2m+3}-2$. 
\end{lemma}

\begin{proof}[Proof of Lemma~\ref{lem:switch}]
We begin by proving that, for all $i\in \Nz$, $x_i=G_i(n)-F_{2m+3}-2$. Indeed, we may update accordingly, for $i\ge 1$: $x_i = x_{i-1} + F_{2n+p_i}$, where $p_i=3$ if $w_i=b$, and  $p_i=4$ if $w_i=a$, where $w_i$ denotes the $i$th entry of the Wythoff word $W=babaababaabaababaababa\ldots$.  This follows, because $G_i(n) = \lfloor i\phi\rfloor F_{2n+2} + iF_{2n+1} + F_{2n+3} -2$, by choosing $a=F_{2n+3}$ and $b=F_{2n+4}$. Hence it suffices to show that (\ref{x}) is equivalent to this update.

In case $j>0$, then, since $j<F_{k-1}$ we get $p_{F_k+j}=x_{F_k+j+1}-x_{F_k+j}=x_{j+1}+F_{2n+k+2}-(x_j+F_{2n+k+2})=x_{j+1}-x_j=p_j$, which is true, by the Fibonacci morphism.

In case $j=0$, then $x_i-x_{i-1} = x_{F_k+1}-x_{F_k} = (x_{F_{k-2}+1}+F_{2n+k+1}) - (x_{F_{k-2}}+F_{2n+k+1})$, since $F_k=F_{k-1}+1+F_{k-2}-1$, and because $x_{F_{k-2}+1}+F_{2n+k+1}=x_{F_k+1}+F_{2n+k+2}$, which holds by plugging in $j=F_{k-1}-1$ in (\ref{x}) and since the updates of the Fibonacci morphism are the same when the parity is the same. We are done with this part of the proof. 

We wish to prove that, for all $k>3$, with $j = 0, \ldots , F_{k-1}-1$, we get $x_{j+1} < F_{2n+k+2}$. The base case $k=3$ gives $j=0$ and therefore $x_{3}=x_{1}+F_{2n+5}$, and clearly $x_1=2F_{2n+3}-F_{2m+3}<F_{2n+5}$. For the general case $k\ge 3$, note that, by letting $j=F_{k-1}-1$ in (\ref{x}), 
\begin{align}\label{xfib}
x_{F_{k+1}}=x_{F_{k-1}}+F_{2n+k+2}. 
\end{align}
Hence it suffices to prove that $x_{F_{k-1}}<F_{2n+k+2}$. Suppose, by induction, that $x_{F_{k-2}}<F_{2n+k+1}$ and $x_{F_{k-3}+1}<F_{2n+k}$. By $F_{2n+k+2}=F_{2n+k+1}+F_{2n+k}$, it thus suffices to prove that , for all $k\ge 3$,
\begin{align}\label{xind}
x_{F_{k+1}}<x_{F_{k}}+x_{F_{k-1}}. 
\end{align}
The base case is that $x_1<x_1+x_0$.  By the induction hypothesis $x_{F_{k-1}}<x_{F_{k-2}}+x_{F_{k-3}}$, we get 
$$x_{F_{k+1}}=x_{F_{k-1}}+F_{2n+k+2}<x_{F_{k-2}}+F_{2n+k+1}+x_{F_{k-3}}+F_{2n+k}=x_{F_{k}}+x_{F_{k-1}},$$ using (\ref{xfib}). 
\end{proof}

Returning to the proof of Theorem~\ref{thm:Q-}, It follows that
\begin{align*}
&F_{2n+2}\lfloor i\phi\rfloor + F_{2n+1}i + F_{2n+3} -2-(F_{2m+3}-2)\in A 
\end{align*}
for all $i, n$ and all $m<n$, so Right cannot find a move from a heap of size $G_i(n)$ to one of size $F_{2m+3}-2$.

For Left, by Theorem~\ref{thm:rcf1} and the induction hypothesis, it suffices to find a move to a heap size of the form $AB(i)+1$, some $i<n$. Hence, we wish to show that $AB(n)-AB(i)-1\in A$, for some $i$. By Corollary~\ref{cor:ab_bb},  we get, for $m>0$ and $0\leq k< F_{2m+1}$, $AB(k+F_{2m})- AB(k)-1=F_{2m+3}-1$ which by Lemma~\ref{lem:aF2n} is equal to $A(F_{2m+2})$. Hence we can take $i = n-F_{2m}$, where $m$ is the largest index such that $n-F_{2m} > 0$.

\end{proof}

In fact, computations indicate that many reduced switches $G_1(n)$ are in canonical form  (even stronger we conjecture ``if and only if''). For example, from a heap of size $8$, Left can move to $4$ and $\heap{4} = \{ 1 |\{ 1 | 0\} \} =_I 1$. 
We may replace Left's option to $\heap{4}$ by 1. We conjecture that analogous replacements holds in general, because Left can move to $\heap{F_{2n+3}-1}$ from $\heap{2F_{2n+3}-2}$. We believe that the Left option, $\heap{F_{2n+3}-1}$, can be replaced by $1$ in general, and all other Left options are dominated or reversed out.

We conjecture that some of the switches are in canonical form.

\begin{conjecture}\label{conj:Q-}
For all $n\in \N$, $\heap{2F_{2n+3} -2} = \cgbgame{1}{s(n)}$.
\end{conjecture}

The motivation for this conjecture is quite strong, only a small part appears to be missing in the proof, so let us sketch some details.
By Lemma \ref{lem:aF2n}, $F_{2n+3} \in B$, so that Right can move to $s(n)$. Since $s(n)$ is increasing, and any other number is 
$ > 2/3$, it suffices to prove that Right cannot move to  $s(m)$ for any $m<n$. Note that, by Corollary \ref{corollary:fibadd}, 
the difference 
\begin{align*}
x=x(m,n) &:=F_{2n+3}-2 - (F_{2m+3}-2)\\ &= F_{2n+3} - F_{2m+3}\\ &=\sum_{m+2\le j\le n+1}F_{2j},
\end{align*}
 is in $A$, for all $0\le m<n$. By Theorem~\ref{thm:RA}, it now follows that $x + F_{2n+3}\in A$, since $x\le F_{2n+3}+1$.
This shows that Right's options are correct. 

For Left, we know that there is a move from $\heap{2F_{2n+3} -2}$ to $\heap{F_{2n+3} -1}$ and $\rcf{F_{2n+3} -1}=1$. We conjecture that this move is reversible. Via the Fibonacci word, one can see that each Right option from $\heap{F_{2n+3}-1}$ is in $B$, and thus each Right option has reduced canonical form $\cgbgame{1}{0}$. In particular Right can move to $\heap{B(1)}=\heap{2}$ and $\heap{2}=\cgbgame{1}{0}$. Hence it is feasible that the Left option $\heap{F_{2n+3}-1}$ reduces to $1$. We know that the reduced canonical form is $1$, so if the infinitesimal $\heap{h^L} - 1$ is nonnegative, for all Left options of $h$, then the conjecture is true. (If it is negative, then Right prefers this option instead.) 

\section{Number theory and combinatorics on words}\label{sec:fibW}

The results, new and old, in this section are presented primarily for use in Section \ref{sec:main}. 

Let $n\in \Nz$. We have that 
\begin{align}\label{bge2}
2\leqslant B(n+1)-B(n)=A(n+1)-A(n)+1\leqslant 3. 
\end{align}

\subsection{Fibonacci representations}
In the sequel we consider positive integers as sums of Fibonacci numbers. More precisely, we represent a positive integer, say $h$, by a multiset of positive integers each of which is taken to be an index of a Fibonacci number. In particular we distinguish $F_1 = 1$ and $F_2=1$. We present three sets of constraints on the multiset of indices such that each set of constraints gives a unique multiset for every positive integers. We allow indices to be unused, used once, or used twice (in which case we call the index a 2). In the first two representations we represent our multisets in binary, and then in the third, in ternary coding.

\noindent{\bf Fact 1}: in the \emph{Zeckendorf representation} of a positive integer $n$, denoted $Z(n)$, no two indices are consecutive and the least index is at least $2$. We may determine the summands for $n$ as follows; the largest Fibonacci number less than $n$, say $f$, is a summand and the remaining summands are the summands in the representation of $n-f$.\\

\noindent{\bf Fact 2}: in the \emph{least-odd representation} of a positive integer $n$, denoted $LO(n)$, indices are distinct, the least index is odd and no two indices are consecutive.\\ 

\noindent{\bf Fact 3}: in the \emph{even representation} of a positive integer $n$, denoted $E(n)$, indices are even, we allow 2s if between each pair of 2s there is at least one unused even-index, $E(n)$. Hence, in ternary coding for example $1020001020=2F_2+F_4+2F_8+F_{10}$ is OK, but neither $2020=2F_2+2F_4 (= 8)$ nor $1020101020 = 2F_2+F_4+F_6+2F_8+F_{10} (= 2+3+8+2\cdot 21+55 = 110)$ is OK. In fact the latter two are uniquely represented as $100000 = F_6$ and $2000000000 = 2F_{10}$. This representation can be obtained by recursively subtracting the largest even-indexed Fibonacci number less than the given number. Note that, for $n>0$ and $k\geqslant 0 $ integers, \[F_{2n} + \left(\sum_{i=0}^k F_{2n+2i}\right) + F_{2n+2k} = F_{2n + 2k +2} + F_{2n-2},\] which explains the uniqueness of the even representation.\\

For example $117 = 89 + 21 + 5 + 2 = F_{11} + F_8 + F_5 + F_3$ in the Zeckendorf representation. Here the least-odd representation coincides with Zeckendorf. In the even representation we rather get $117=55+55+3+3+1 = 2F_{10} + 2F_4 + F_2$. On the other hand, for example  $11=F_6 + F_4 = F_6+F_3+F_1$. So here Zeckendorf coincides with the even representation, whereas the least-odd representation differs.

We use the following folklore result, which has appeared in several versions since the 70s, but one nice source is \cite{S1}.

\begin{proposition}\label{prop:zeck}
The unique Zeckendorf representation of a positive integer $x$, $Z(x)$, ends in an even number of 0s if and only if $x\in A$. Precisely, $n\in \N$ is the right shift of $A(n)$, whereas $B(n)$ is the left shift of $A(n)$. Unless $Z(A(n))$ ends in zero 0s the representation of $n$ is $Z(n)$; otherwise it is the least-odd representation. 
\end{proposition}
This classical result is interesting in several ways. Indeed, one purpose in our setting is to find, in linear time, whether a given number is a legal move or not. Obviously we do a standard Zeckendorf decomposition of the given number, and read off the last few digits, that is if the least Fibonacci number is even-indexed, then the number belongs to $A$, and otherwise to $B$. That is, if there is an even number of right most 0s, of a given number in the Zeckendorff representation, it corresponds to a legal move for Left, otherwise for Right. In this paper,  it will also be indispensable to quickly find the correct type of strategy. This is where we can use that $B(n)$ is the left shift of $A(n)$ and $n$ is its  right shift. 
For example if $A(n)=9=10001_{\text{Zeck}}$ then $n=6=10001_{\text{LO}}$, where the last $1=F_1$, so that $n$ is written in the least-odd representation. We will be interested in numbers of the form $AB(n)+1$. To find out whether $x$ is of this form, we thus investigate $Z(x-1)$. If its right shift ends in an odd number of 0s, then $x=AB(n)+1$ for some $n$. That is, $x$ is of the correct form if and only if $Z(x-1)$ ends in a strictly positive even number of 0s. We will sometimes use the notation $z_1(n)$ for the least index of a Fibonacci number in the Zeckendorf representation of $n$.

Given $x\in \N$, let us sketch an algorithm for transforming $Z(x)$ to the even representation, $E(x)$. Let $E_0(x) := Z(x)$. At step $i\geqslant 0$, denote the index of the least odd-indexed Fibonacci number in $E_i(x)$ by $n = n(i)$. If there is none, we are done; take output $E(x)=E_i(x)$. If $n=3$ then put $F_3\rightarrow 2F_2$ and let the rest of the representation be as in $E_i(x)$; return $E_{i+1}(x)$. If $n>3$ and $F_{n-3}$ is also present put $F_n+F_{n-3}\rightarrow 2F_{n-1}$, otherwise put $F_n \rightarrow F_{n-1} + F_{n-2}$, and in either case let the rest of the representation be as in $E_i(x)$; return $E_{i+1}(x)$. We call this algorithm, the ZE-algorithm. It is clear that it terminates with the unique even representation of $x$. We will also need the reverse to this algorithm in the section for the numbers, but the application is fairly special and best presented in that context; we will also use the following lemma.

\begin{lemma}\label{evenreprZodd}
Let $n\in \N$. Then $z_1(n)$ is odd if and only if the smallest $n$ terms in the even representation of $n$ are
\begin{align}\label{z1even}
F_2+F_4+\cdots +F_{2i} + 2F_{2i+2}, \text{ for some } i\ge 0.
\end{align} 
\end{lemma}
\begin{proof}
Suppose that (\ref{z1even}) holds. Then if $i>0$, by successively applying the formula 
\begin{align}\label{doubltransf}
2F_{2i+2}=F_{2i}+F_{2i+3}, 
\end{align}
the claim holds, because $2F_2=F_3$.  Note that the greater terms will not affect (the parity of) the smallest index.

Suppose that (\ref{z1even}) does not hold. Then, note that unless there is any 2, $z_1(n)$ is even because the representations coincide. If there is a 2, by successively applying (\ref{doubltransf}), by the even representation (there is a gap somewhere between each pair of 2s) it will terminate with a single smallest even Fibonacci number (and perhaps there are smaller even Fibonacci terms but they will then remain unaffected by the transformation of the greater terms). 
\end{proof}

In Section \ref{sec:numbers} we will use this lemma.

\begin{lemma}\label{lem:repr}
Each number in $\widehat{B^2}$ has a `4' as the least index in the even representation.
\end{lemma}
\begin{proof}
Write $n$ in the least-odd representation and left shift 4 times. We get $B^2(n) $ in the Zeckendorf representation as $Z(n)000$, that is the least Fibonacci term has odd index $\geqslant 5$. That is, we get $B^2(n) + 1 = Z(n)001$ with the second least Fibonacci term $F_{2k+1}\geqslant F_5$. In case $k=2$ then  $F_{5}+1 = 2F_{4}$, otherwise, if $k>2$, then $F_{2k+1}+1=F_{2k}+F_{2k-2}+\ldots +F_4+F_3+1=F_{2k}+F_{2k-2}+\ldots +F_4+F_4$. We need to prove that the larger terms cannot affect the occurrence of (the right most) $F_4$.
By the ZE-algorithm, in $Z(B^2(n)+1)\rightarrow E(B^2(n)+1)$, at most one of the $F_4$'s can be converted (to $F_6$) and there will be no introduction of $F_2$ since the possible odd-indexed terms are $\geqslant F_5$. By the uniqueness of the even representation, we are done.
\end{proof}

\subsection{The Fibonacci morphism}
We use \emph{words} and prove statements about them. We use \cite{L} as a source of definitions for \emph{words}.

The Fibonacci sequence was defined in the introduction. The \emph{Fibonacci morphism} $\varphi:\{a,b\}^\star \rightarrow \{a,b\}^\star$ is 
\[\varphi: \left\{
    \begin{array}[l]{l}
      a\rightarrow ab\\
      b\rightarrow a
    \end{array}\right.
\]
Here $\varphi^n(x):=\varphi(\varphi^{n-1}(x))$, for $n>0$ and all words $x\in \{a,b\}^\star$, also $\varphi^0(x) := x$. We are in particular interested in the following recurrence:
  \begin{align}
    \varphi^0(a)&=a\label{fig:morphism}\\
    \varphi^1(a)&=ab\notag\\
    \varphi^2(a)&=aba\notag\\
    \varphi^3(a)&=abaab\notag\\
    \varphi^4(a)&=abaababa\notag\\
    \varphi^5(a)&=abaababaabaab\notag\\
    \varphi^6(a)&=abaababaabaababaababa\notag\\ \notag
  \end{align}
and so on.

For all $n \geqslant 2$,
\begin{equation}\label{n1n2}
\varphi^n(a) = \varphi^{n-1}(a)\varphi^{n-2}(a)
\end{equation}
from which we conclude that a prefix of $\varphi^n(a)$ is a prefix of $\varphi^m(a)$ for all $n \leqslant m$. We can thus define the (infinite) Fibonacci word:
\[\varphi^\infty := \lim_{n \to \infty}\varphi^{n}(a)=abaababaabaababaababa\ldots\]

For readability, we sometimes write $\varphi^n$ instead of $\varphi^n(a)$. It is handy to have
the following variation of (\ref{n1n2}). For all $n \geqslant 3$,
\begin{align}\label{varphi:form}
\varphi^{n} = \varphi^{n-2}\varphi^{n-3}\varphi^{n-2}.
\end{align}
For a finite word $x$ we use $|x|$ to denote the length (i.e.~the total number of letters) in the word $x$.
For a finite word $x$ we use $|x|_a$ to denote the number of occurences of $a$ in the word $x$.
For any word $x$ we let $x_n$ denote the prefix of $x$ of length $n$ (if it exists). Later, in our application, we will often use the \emph{Wythoff word} $W:=b\varphi^\infty$.

\begin{lemma}\label{lem:varphicount}
For a finite word $x$ on the alphabet $\{a,b\}$,
\begin{align*}
|\varphi(x)|_b &= |x|_a\\
|\varphi(x)|_a &= |x|.
\end{align*}
\end{lemma}
\begin{proof}
From the Fibonacci morphism: $|\varphi(x)|_b = |x|_a$; $|\varphi(x)|_a = |x|_a + |x|_b$, which is $|x|$ as our word is on the alphabet $\{a,b\}$.
\end{proof}

\begin{lemma}\label{lem:fibprecounts}
For $n \geqslant 0$,
\begin{align*}
|\varphi^n|_{b} &= F_{n}\\
|\varphi^n|_{a} &= F_{n+1}\\
|\varphi^n|_{\phantom{2}} &= F_{n+2}.\\
\end{align*}
\end{lemma}

\begin{proof}
This is true for $n = 0$. Using Lemma \ref{lem:varphicount} and the Fibonacci identity we get the result by induction.
\end{proof}

\begin{proposition}\label{varphi:double}
If $n\geqslant 2$ then \[{\varphi^\infty}_{2F_{n+2}}=\varphi^n\varphi^n.\] 
\end{proposition}

\begin{proof}
Suppose that $n=2$. Then $\varphi^4=\varphi^2\varphi^b\varphi^2=abaababa=\varphi^2\varphi^2ba$, which proves this case since $|\varphi^2\varphi^2|=6=2F_4$.
In case $n\ge 3$, we apply (\ref{varphi:form}) twice and then (\ref{n1n2}),
\begin{align*}
\varphi^{n+2} 
&= \varphi^{n}\varphi^{n-1}\varphi^{n}\\
&= \varphi^{n}\varphi^{n-1}\varphi^{n-2}\varphi^{n-3}\varphi^{n-2}\\
&= \varphi^{n}\varphi^{n}\varphi^{n-3}\varphi^{n-2}.
\end{align*}
By Lemma~\ref{lem:fibprecounts}, the length of the prefix is correct.
\end{proof}
We call the following lemma the \emph{glueing principle}. 

\begin{lemma}\label{lem:obvious}
Let $x$, $x'$, $y$, and $z$ be finite words. Suppose that $y$ is a prefix of $x$ and $z$ is a suffix of $x'$. If $x = x'$ and $| y | + | z | =| x |$ then $x = yz$.
\end{lemma}
\begin{proof}
Obvious.
\end{proof}
The way we will use it is as follows. We will find words of the form $xyx$. Then Any factor of this word of length $|xy|$ contains exactly the same number of letters of each kind.
\begin{definition}
For non-negative integers $n$ and $m$ and a finite word $x$ on the alphabet $\{a,b\}$, let \[S_{n,m}(x) = n|x|_b +m|x|_a.\]   
\end{definition}

\begin{lemma}\label{lem:S12}
 \[ S_{1,2}(\varphi^n) = F_{n+3}.\]
\end{lemma}
\begin{proof}
Using Lemma~\ref{lem:fibprecounts}, $S_{1,2}(\varphi^n) =  1|\varphi^n|_b +2|\varphi^n|_a = F_n + 2F_{n+1}  =  F_{n+3}.$
\end{proof}

\begin{theorem}\label{thm:RA}
If $k\leqslant F_{n+3} + 1$ and $n\geqslant 2$ then $k+F_{n+3} \in A$ if and only if $k\in A$.
\end{theorem}

\begin{proof}
By Lemma \ref{lem:S12}, $S_{1,2}(\varphi^n)=F_{n+3}$ for $n\geqslant 0$.
By Proposition~\ref{varphi:double}, for all $n\geqslant 2$, each prefix of $b\varphi^n\varphi^n$ is a prefix of $W$. 
Hence, for $n\geqslant 2$, by Lemma \ref{lem:obvious}, if $1\leqslant k\leqslant S_{1,2}(b\varphi^n)$, then $k+S_{1,2}(\varphi^n)\in A$ if and only if $k\in A$. 
\end{proof}

\begin{proposition}\label{prop:Aif0}
If $t$ is the $k$th letter of $\varphi^{\infty}$, then
$\begin{cases}
k \in A, & \text{if } t = a\\
k \in B, & \text{if } t = b.\\
\end{cases}$
\end{proposition}

\begin{proof}
  We prove the result by showing it holds for every prefix, that is, it holds for $\varphi^n$ for every $n$. It holds for $n = 1$; we prove the rest by induction.
 Every $1$ in $\varphi^{k+1}$ comes from an $a$ in $\varphi^{k}$ via the $a \to ab$ part of the morphism; the $b$ is not simply one space to the right of where the preimage $a$ was because each of the $i-1$ preceding $a$s in the pre-image became $ab$. Hence, the location of the $i$th b in $\varphi^{k+1}$ is the location of the $i$th $a$ in $\varphi^{k}$ plus 1 plus $i-1$, which is $A(i) + i$ by induction. Furthermore $A(i) + i = B(i)$. The location of $a$s is given by $A$ using complementarity.
\end{proof}

\begin{corollary}\label{cor:notremove}
Let $p$ be a prefix of $\varphi^\infty$.
\begin{itemize}
\item If $p$ ends in $a$ and $|p|_a = n$ then $|p| = A(n)$.
\item If $p$ ends in $b$ and $|p|_b = n$ then $|p| = B(n)$.
\end{itemize}
\end{corollary}

\begin{proof}
This follows from Proposition~\ref{prop:Aif0} because both of the sequences $A$ and $B$ are increasing.
\end{proof}

\begin{lemma}\label{lem:Sphi}
If $x$ be a finite word on the alphabet $\{a,b\}$ then $S_{F_{i},F_{i+1}}(\varphi(x)) = S_{F_{i+1},F_{i+2}}(x)$.
\end{lemma}

\begin{proof}
Using Lemma~\ref{lem:varphicount},
\begin{equation*}
\begin{aligned}S_{F_{i},F_{i+1}}(\varphi(x)) 
&= F_{i}|\varphi(x)|_1 +F_{i+1}|\varphi(x)|_0\\
&= F_{i}|x|_0 +F_{i+1}(|x|_0+|x|_1)\\
&= F_{i+2}|x|_0 +F_{i+1}|x|_1\\
&= F_{i+1}|x|_1 + F_{i+2}|x|_0\\
&= S_{F_{i+1},F_{i+2}}(x)\\
\end{aligned}
\end{equation*}
\end{proof}

\begin{corollary}\label{cor:16}
$|\varphi(x)| = S_{1,2}(x)$
\end{corollary}

\begin{proof}
Using Lemma \ref{lem:Sphi} and letting $i = 1$ we get,
$S_{1,1}(\varphi(x)) = S_{1,2}(x)$
\end{proof}

\begin{lemma}\label{lem:17}
For all $n$,  $S_{1,2}(W_{n}) = A(n)$ and $S_{2,3}(W_{n}) = B(n)$.
\end{lemma}

\begin{proof}
  The first part is well known and the second part is easily derived from noting that both integer inputs to $S$ are exactly 1 more, leading to counting 1 extra each of $n$ times and (\ref{bann}).
\end{proof}

It is convenient to regard the sequences $A$ and $B$ as functions on the non-negative integers. We write $AB(k)$ rather than $A(B(k))$, $A^n(k)$ rather than $A(A^{n-1}(k))$, and so on. In fact, we will later use a word notation for this composition of functions. It will be clear from the context, whether we discuss Sturmian words (as previously in this section) or whether we use the word notation for composition. In Table~\ref{tab:sequences2} we show our two most relevant compositions of the $A$ and $B$ sequences.

\begin{table}[ht]
\begin{center}
\begin{tabular}{|l||c|c|c|c|c|c|c|c|c|c|c|c|c|c|c|}\hline
$n$      & 0& 1& 2& 3& 4& 5& 6& 7& 8& 9&10&11&12&13&14\\ \hline\hline
$A(n)$    & 0& 1& 3& 4& 6& 8& 9&11&12&14&16&17&19&21&22\\ \hline
$B(n)$    & 0& 2& 5& 7&10&13&15&18&20&23&26&28&31&34&36\\ \hline
$AB(n)$& 0& 3& 8&11&16&21&24&29&32&37&42&45&50&55&58\\ \hline
$B^2(n)$& 0& 5& 13&18&26&34&39&47&52&60&68&73&81&89&94\\ \hline
$W$    &$b$& $a$& $b$& $a$& $a$& $b$& $a$& $b$& $a$& $a$& $b$& $a$& $a$& $b$& $a$\\ \hline
\end{tabular}\caption{The table displays the first few terms of the $A$ and $B$ sequences together 
with the corresponding prefix of the infinite Wythoff word $W=b\varphi^\infty$ (obtained via the Fibonacci word $\varphi^\infty$), defined in Section \ref{sec:fibW}. The sequences $AB(n)=A(n)+B(n)$ and $B^2(n)=B(B(n))$ are of special interest. A calibration of the latter three sequences gives the four sets that partition the nonnegative integers in Table~\ref{tab:partitioning}.}\label{tab:sequences2}
\end{center}
\end{table}

We have the following results from Kimberling, \cite{K}. For all $n$,
\begin{align}
  A^2(n) &= B(n) - 1 \label{Ki1}\\
  AB(n) &= A(n) + B(n) \label{Ki2}\\
  BA(n) &= A(n) + B(n) -1\label{Ki3}\\
  B^2(n) &= A(n) + 2B(n).\label{Ki4}
\end{align}

The Fibonacci numbers alternate between being in $A$ and $B$.
For example, $F_2=A(1)=1$, $F_3 = B(1) = 2$, $F_4=AB(1)=3$, $F_5=B^2(1)=5$, and in general:

\begin{lemma}\label{lem:aF2n}
\begin{align}\label{FfromF2n-1}
A(F_{2n-1})&=F_{2n}, \\
B(F_{2n-1})&= F_{2n+1} ,\\
A(F_{2n}) &= F_{2n+1}-1,\\ 
B(F_{2n}) &= F_{2n+2}-1. \label{BF2n}
\end{align}
\end{lemma}
\begin{proof}
Combine Corollary \ref{cor:16} with Lemma \ref{lem:17}.
\end{proof}

We can also write the previous results in long form:
\begin{align}\label{F2nF2n1}
F_{2n} &= AB^{n-1}(1)\mbox{ and,}\\
F_{2n+1} &= B^n(1).
\end{align}

Kimberling generalizes the Equations (\ref{Ki1}) to (\ref{Ki4}) and in particular proves (a more general form of) the following theorem on complementary Wythoff sequences.

\begin{theorem}(\cite[Theorem 1 and Corollary 2]{K})\label{theorem:genK}
For all $n\geqslant 1$ and $i\geqslant 0$, 
\begin{align*}
F_{2n-3}A(i) + F_{2n-2}B(i) &= AB^{n-1}(i),\\
F_{2n-2}A(i) + F_{2n-1}B(i) &= F_{2n}A(i) + F_{2n-1}i\\ 
&= B^n(i),
\text{ and }\\
F_{2n}A(i) + F_{2n+1}B(i)-F_{2n+1} &= AB^{n}A(i).\\
\end{align*}
 
\end{theorem}

This result implies that, for all $n\geqslant 0$ and $i>0$, 
\begin{align}\label{ABnABn}
AB^n(i)-AB^n(i-1)\in\{F_{2n+2}, F_{2n+3}\}
\end{align}
and if $n>0$
\begin{align}\label{BnBn}
B^n(i)-B^n(i-1)\in\{ F_{2n+1}, F_{2n+2}\}.
\end{align}
Precisely 
\begin{align}\label{ABnABn2}
AB^n(A(i))-AB^n(A(i)-1)&=F_{2n+2},\\
AB^n(B(i))-AB^n(B(i)-1)&=F_{2n+3},
\end{align}
and 
\begin{align}\label{BnBn2}
B^n(A(i))-B^n(A(i)-1) &= F_{2n+1},\\
B^n(B(i))-B^n(B(i)-1) &= F_{2n+2}.
\end{align}

Clearly such assertions can be generalized, see \cite{K}. Any permutation of the letters $A$s and $B$s (regarded as composition of functions) produces the same pair of \emph{gaps} between consecutive entries in the corresponding sequence,  only the initial entry may differ. For example, we also have, for all $n\geqslant 0$ and $i>0$, 
\begin{align}\label{BnABnA}
B^nA(A(i))-B^nA(A(i)-1)&=F_{2n+2},\\
B^nA(B(i))-B^nA(B(i)-1)&=F_{2n+3}.
\end{align}
 
The next result generalizes the above lemmas in another direction. We will only use a corollary of this result, stated just after its proof.

\begin{theorem}\label{lemma:fibadd}
Let $g_c:=\sum c_iF_{2i}$  and $h_c:=\sum c_iF_{2i-1}$ where $c_i\in\{0,1,2\}$, $i\in \mathbb{Z}_{> 0}$, $c=(c_i)$, $c_i=0$ for all but finitely many $i$, and where precisely one of the following holds
\begin{itemize}
\item[(i)] $c_i \in \{0,1\}$ for all $i$ (in particular implying no repetition of Fibonacci numbers); 
\item[(ii)] $c_1=0$ (otherwise no restriction). 
\end{itemize} 
Then $g_c=A(h_c)$, that is  
\begin{align}\label{aF}
\sum_i c_iA(F_{2i-1}) = A\left(\sum_i c_iF_{2i-1}\right). 
\end{align}
\end{theorem}

\begin{proof}
We begin by noting that, by Lemma \ref{lem:aF2n}, $g_c = A(h_c)$ implies (\ref{aF}).

For (i), let $g(k)=g_c(k)$ denote a sum of even-indexed distinct Fibonacci numbers each with index less than or equal to $2k$. Then $g(k)$ is strictly less than $F_{2k+1}$ with the corresponding $h(k)=h_c(g(k))<F_{2k}$. That is, given  $g(k)$, we also know $h(k)$. The induction hypothesis is that, for all $g(k)$, for a given $k$,
\begin{align}\label{gkSp}
g(k)=S_{1,2}W_{h(k)}. 
\end{align}
That is, by $A(h(k)) = S_{1,2}W_{h(k)}$, we assume that, for this given $k$, the additive property (\ref{aF}) for the $A$ sequence will be  satisfied. Note that $W_{h(k)}$ will be the prefix of $b\varphi^{2k-2}$ of length $h(k)$. Hence, we want to show that any $g(k+1)$ is a partial  sum of $$S_{1,2}(b\varphi^{2k}) = S_{1,2}(b\varphi^{2k-2}\varphi^{2k-3}\varphi^{2k-2}),$$ with the prescribed $h(k+1)$ initial terms. We may assume that $g(k+1)$ has $F_{2k+2}=A(F_{2k+1})$ as a term, since otherwise we are done (because it would also satisfy the condition for $g(k)$). This means that $F_{2k+1}$ is a term in $h(k+1)$. The factor $\varphi^{2k-1}=\varphi^{2k-2}\varphi^{2k-3}$ has length $F_{2k+1} = F_{2k} + F_{2k-1}$. 
Now the induction hypothesis gives the result, since $W_{h(k)}$, without the initial $b$, is a prefix, say $p$, of the suffix $\varphi^{2k-2}$ of $\varphi^{2k}$. That is, by (\ref{gkSp}), we get $g(k+1) = S_{1,2}\varphi^{2k-1} + g(k) = S_{1,2}(b\varphi^{2k-1}p)$, where $b\varphi^{2k-1}p$ is the prefix of $b\varphi^{2k}$ of length $h(k+1)=F_{2k+1}+h(k)$. This ends the first part of the proof.
 
A similar idea gives the proof for case (ii). Given $c$, let $g(k)$ denote a sum of even-indexed Fibonacci numbers of the form in (ii), each with index less than or equal to $2k$. Then $g(k)$ is less than $2F_{2k+1} = 2F_{2k} + 2F_{2k-1}$ with the corresponding $h(k)$ less than $2F_{2k}$. The inductive hypothesis is that $g(k)=S_{1,2}p$, where $p$ is a prefix of $b\varphi^{2k-2}\varphi^{2k-2}$ of length $h(k)$, that is, the additive hypothesis (\ref{aF}) for the $A$ sequence is satisfied.  Hence, we want to show that $g(k+1)$ is a partial  sum of $$S_{1,2}(b\varphi^{2k}\varphi^{2k})$$ with the prescribed $h(k+1)$ terms. We may assume that $g(k+1)$ has $S_{1,2}(\varphi^{2k-1}) = F_{2k+2}=A(F_{2k+1})$ as a term, since otherwise we are done. Note that 
\begin{align}\label{vvvvvvvv}
\varphi^{2k}\varphi^{2k}&=\varphi^{2k-2}\varphi^{2k-3}\varphi^{2k-2}\varphi^{2k-2}\varphi^{2k-3}\varphi^{2k-2}.
\end{align}
Hence, if the prefix $p$ of $b\varphi^{2k-2}\varphi^{2k-2}$ is not a prefix of $b\varphi^{2k-2}$, then the argument goes through irrespective of whether we add $F_{2k+2}$ or $2F_{2k+2}$, namely the second last or the last $\varphi^{2k-2}$ in (\ref{vvvvvvvv}) will fill the role of the last one from the hypothesis. Thus, the remaining case is if $p$ is  a prefix of $b\varphi^{2k-2}$, but we add $2F_{2k+2}$. The problem is that $\varphi^{2k-2}=\varphi^{2k-3}\varphi^{2k-4}\ne \varphi^{2k-4}\varphi^{2k-3}$. However if we erase the two letter suffixes, that is $ba$ and $ab$ respectively, then the remaining prefixes are the same (since they are palindromes). Hence the only problem is for $A(F_{2k})$ since $A(F_{2k})+2=F_{2k+1}+1=S_{1,2}(b\varphi^{2k-2})$. However, $F_{2k+1}-1$ is an illegal configuration in the setting of (ii) since
it is not a sum of even-indexed Fibonacci numbers all $\geqslant 3$.
\end{proof}

Note that, given (i) or (ii), it follows, by (\ref{aF}) and by $B(n)=A(n)+n$, for all $n$, that 
\begin{align*}
\sum_i B(c_iF_{2i-1}) = B\left(\sum_i c_iF_{2i-1}\right).
\end{align*}
Note that the identities for $g_c$ and $h_c$ in the theorem only hold given certain conditions on $c$, specifying special Fibonacci representations of certain integers. Below we also consider something called the unique ``even Fibonacci representation" of any positive integer, but this is an entirely different story.

\begin{corollary}\label{corollary:fibadd} 
For all $0<m<n$, 
 $F_{2m} + F_{2m+2} + \ldots + F_{2n}= A(F_{2m-1})+A(F_{2m+1})\ldots + A(F_{2n-1})\in A$ and  if in addition $1<m$ then $F_{2m} +F_{2m+2}+ \ldots + 2F_{2n} = A(F_{2m-1}) + A(F_{2m+1})\ldots + 2A(F_{2n-1})\in A$.
\end{corollary}
\begin{proof}
Apply Theorem \ref{lemma:fibadd}, (\ref{aF}).
\end{proof}
The following lemma regards a certain palindrome structure of the Fibonacci word. A proof can be found in \cite[Section 2.2]{L}, in particular their Example~2.2.2 together with Proposition~2.2.4.
\begin{lemma}\label{lem:palind}
Any prefix of the infinite Fibonacci word $\varphi^\infty$ of length $F_n-2$ is a palindrome.
\end{lemma}
We also reframe Lemma \ref{lem:palind} in the setting of \WPS.
\begin{corollary}\label{cor:palind}
If $w = W_{F_n-1}$ then $wb$ is a palindrome.
\end{corollary}

If $n$ is even, then $W_{F_n}$ is a palindrome. This follows since the even recurrences of $\varphi$ end in $ba$ (by induction using that the odd recurrences end in $ab$). 

Let $w = W_{F_n-1}$ for some $n$.

The \emph{palindrome principle} is the following.
Let $p$ be a palindrome. If a player (for example Right) moves from a heap of size $S_{2,3}(wb)$, then he moves to a heap with size in $B$; his moves demonstrates a split $p = xy$ such that $S_{2,3}(x) \in B$ iff $S_{2,3}(y) \in B$.

The next results concern an \emph{invariance principle} of the Fibonacci word.

\begin{lemma}\label{lem:prefixfactor}
If $n \geqslant 2$ and $v$ is a factor of $\varphi^{n}$ with $|v| = F_{n+1}$ then $|v|_b = F_{n-1}$.
\end{lemma}

\begin{proof}
Recall  (\ref{varphi:form}) and Lemma \ref{lem:obvious}. The length of the factor forces it to contain all of the $\varphi^{n-3}$ factor and the rest of our factor is some sufffix of $\varphi^{n-2}$ on the left and the remainder of $\varphi^{n-2}$ as a prefix on the right. That is, any factor of length $F_{n+1}$ contains $F_{n-3} + F_{n-2} = F_{n-1}$ $b$s.
\end{proof}

What we want is the following consequence:

\begin{corollary}\label{cor:prefixfactor}
If $n \geqslant 2$ and $v$ is a factor of $\varphi^\infty_{F_{n+2} - 2}$ with $|v| = F_{n}$ then $|v|_b = F_{n-2}$.
\end{corollary}

\begin{proof}
We combine Lemmas~\ref{lem:palind}~and~\ref{lem:prefixfactor}. Let $w = \varphi^\infty_{F_{n+2} - 2}$. Each factor of $w$ of length $F_{n}$ is nearly contained in the prefix $\varphi^{n-1}$ of $w$ of length $F_{n+1}$, or in its suffix of length $F_{n+1}$. In fact, the suffix and prefix share precisely the middle $F_{n-1}+2$ letters. By Lemma \ref{lem:palind}, the factor with $F_{n-2}-2$ letters that follows immediately to the right of $\varphi^{n-2}$ constitutes a palindrome. (Indeed this factor is a prefix of $\varphi^{n-4}$ which in its turn is a prefix of $\varphi^{n-3}$ and $\varphi^{n-1}=\varphi^{n-2}\varphi^{n-3}$.) Hence, by applying Lemma \ref{lem:prefixfactor}, any factor with  length $F_n=F_{n-1}+2+F_{n-2}-2$ has the correct number of $b$s.
\end{proof}

For example, 
$n=2$ gives $w= 0$, $|v|=1$, $|v|_b=0$; 
$n=3$ gives $w=aba$, $|v|=2$, $|v|_b=1$; 
$n=4$ gives $w=abaaba$, $|v|=3$, $|v|_b=1$, and so on. 

The next result can be generalized (by using the setting of Lemma~\ref{lem:32} and using the main theorem of \cite{K}), but we settle with a form that is sufficient for this work.
\begin{corollary}\label{cor:ab_bb}
Let $n>0$. Then 
\begin{itemize}
\item[(i)] for $0<k< F_{2n+1}$, $B^2(k+F_{2n})- B^2(k)=F_{2n+4}$;
\item[(ii)] for $0\leqslant k< F_{2n}$, $B^2(k+F_{2n-1})- B^2(k)=F_{2n+3}$;
\item[(iii)] for $0 < k< F_{2n+1}$, $AB(k+F_{2n})- AB(k)=F_{2n+3}.$
\item[(iv)] for $0\leqslant k< F_{2n}$, $AB(k+F_{2n-1})- AB(k)=F_{2n+3}.$
 \end{itemize}
 
\end{corollary}
\begin{proof} For (i), consider the prefix $\varphi^{2n}$ of $\varphi^\infty$, of length $F_{2n+2}$. We get  $B^2(|\varphi^{2n}|+1)=S_{5,8}(b\varphi^{2n})$. By $W=b\varphi^\infty$, the inequality $k+F_{2n}<F_{2n+2}$, lets us disregard the 2-letter suffix of  $\varphi^{2n}$, in the $B^2$ word. Then, by Corollary~\ref{cor:prefixfactor}, the invariance principle applies and so, each remaining factor of length $F_{2n}$ has the same sum, namely $5F_{2n-2}+8F_{2n-1}=F_{2n+4}$. Here, we use (\ref{BnBn2}). Note that the case $k=0$ is not included in Corollary~\ref{cor:prefixfactor}, so it has to be checked separately using equations in this section. For items (ii) and (iv), the cases $k=0$ can be verified to hold, and the rest  (including case (iii)) is analogous with (i). 
\end{proof}

Let $\{c\} := c-\lfloor c\rfloor$ denote the \emph{fractional part} of a real number $c$. Thus $0\leqslant \{c\}<1$. This is a well known result, but we prove it for completeness.
 
\begin{lemma}\label{lem:frac}
For all positive integers $n$, $\{\phi\lfloor \phi^2 n\rfloor\}<\phi^{-2}$ and $\{\phi\lfloor \phi n\rfloor\}>\phi^{-2}$.
\end{lemma}
\begin{proof} 
The expressions are true for $n=1$. Suppose that they are also correct for all positive integers less than or equal to $m$, where $m=\lfloor \phi^2 n\rfloor$, some $n$. Then in particular $\{\phi\lfloor \phi^2 n\rfloor\}<\phi^{-2}$ and, by complementarity and (\ref{bge2}), there is an $i$ such that $\lfloor \phi n\rfloor + 1 = \lfloor \phi i \rfloor$, and it follows, by a simple modulo 1 estimate, that $\{\phi\lfloor \phi i\rfloor\}>\phi^{-2}$. 

If, on the other hand $m=\lfloor \phi n\rfloor$, some $n$, then there are two possibilities for $m+1$. Either $m+1=\lfloor \phi (n+1)\rfloor$ or $m+1=\lfloor \phi^2 i\rfloor$, some $i<n$. Hence, the second case is immediate by the induction hypothesis. For the first case, we use some algebra on the fractional parts. To begin with note that $\lfloor \phi n\rfloor+1 = \lfloor \phi (n+1)\rfloor$ implies that $\{\phi n\}+\phi^{-1}=\{\phi(n+1)\}$ and therefore also that 
\begin{align}\label{phi-2}
\{\phi n\}\in (0,\phi^{-2}). 
\end{align}
Let $g(n):=\phi^2\{\phi\lfloor \phi n\rfloor\}$. Then, by the hypothesis, $g(n)>1$ and we have to show that $g(n+1)>1$. This is equivalent to showing that 
\begin{align}\label{phi2phi-2}
g(n)\in (\phi^2-\phi^{-1},\phi^2). 
\end{align}
We have that
\begin{align*}
g(n)&=\phi^2\{\phi^{-1}\lfloor \phi n\rfloor\}\\
&=\phi^2\{\phi^{-1}(\phi n-\{\phi n\})\}\\
&=\phi^2-\phi^2\{\phi^{-1}\{\phi n\}\}.
\end{align*}
Then, by (\ref{phi-2}), we get (\ref{phi2phi-2}). 
\end{proof}

Note that, for the fractional part and $r$ a real number, 
\begin{align}\label{fracpart}
\{-r\}=1-\{r\}.
\end{align}
We also have, for real $r$ and $s$, if $0\le\{r\}-\{s\}$ then $\{r\}-\{s\}= \{r-s\}$ and also, by (\ref{fracpart}), 
\begin{align}\label{fracpart2}
0\geqslant \{s\}-\{r\}= \{s-r\}-1.
\end{align}

\begin{lemma}\label{lem:not1} 
Let $h\in \mathcal N$. Then $h-x\not \in A$ if $x\in \widehat{AB_0}$.
\end{lemma}

\begin{proof}
Since $\n = \{F_{2n +3} - 2\mid n\geqslant 0\}\cup\{3\lfloor n\phi\rfloor + 2n+1\mid n\geqslant 0\}$, it suffices to demonstrate that 
\begin{align}\label{bla}
2B(n)+A(n)+1-(B(i)+A(i)+1)=B^2(n)-AB(i)\not\in A 
\end{align}
and 
\begin{align}\label{check1}
F_{2n+1} - 2-(B(i) + A(i) + 1) &= F_{2n+1} - AB(i) - 3 \notag\\ 
&= B(F_{2n-1}) - AB(i) - 3\not\in  A,
\end{align}
for all $n>0, i\geqslant 0$. 

We begin by showing that
\begin{align}\label{blab}
B(B(n)-A(i))=B^2(n)-BA(i)-1, 
\end{align}
whenever 
\begin{align}\label{BA0}
B(n)-A(i) > 0, i>0,
\end{align}
which implies (\ref{bla}). 

Note that, by definition and by $\lfloor \phi^2 x \rfloor=\lfloor \phi x \rfloor+x$ for all nonnegative integers $x$, (\ref{blab}) is equivalent to 
$$\lfloor \phi (B(n)- A(i))\rfloor = \lfloor \phi B(n)\rfloor-\lfloor \phi A(i)\rfloor - 1 ,$$ 
for all $i$ and $n$. 

By Lemma \ref{lem:frac}, we have that 
\begin{align}\label{blaha}
\{\phi B(n)\}<2-\phi <\{\phi A(i)\}, 
\end{align}
for all $n$ and $i$. Let $c>d$ be positive integers. By properties of Beatty sequences we have that $\lfloor \phi (c-d) \rfloor \in \{\lfloor \phi c \rfloor - \lfloor \phi d \rfloor -1, \lfloor \phi c \rfloor - \lfloor \phi d \rfloor \}$. We will show that, with $c=B(n)$ and $d=A(i)$, the first element in this set will be attained, which suffices to prove (\ref{blab}). To this purpose, note that, by definition of integer part, $\lfloor \phi (c-d) \rfloor = \lfloor \phi c \rfloor - \lfloor \phi d \rfloor -1$ is equivalent to $\phi c-\phi d-\{\phi (c-d)\}=\phi(c-d)+\{\phi d\}-\{\phi c\}-1$, which is equivalent to that  $-\{\phi (c-d)\} = \{\phi d\} - \{\phi c\} - 1$. By combining (\ref{fracpart2}), (\ref{BA0}) and (\ref{blaha}), we have that  
\begin{align*}
-\{\phi (B(n)-A(i))\}&= \{\phi A(i)\} - \{\phi B(n)\} - 1, 
\end{align*}
which proves the first part, (\ref{bla}).

For (\ref{check1}), by the first part, it suffices to demonstrate that $B(F_{2n-1}-A(i))-3\in B$, for all $i$ and $n$. This holds if the prefix of $W$ of length $F_{2n-1}-A(i)>1$ ends in a 0, for all $i$ and $n$. Since we begin with an odd-indexed Fibonacci number, the result follows by the palindrome principle and the observation in Corollary~\ref{cor:notremove}, that the number of letters in $\varphi^\infty$ to the left and including the $n$th 0 is in $A$. 
\end{proof}

\subsection{Two consecutive Right moves}
The following result is used, see Theorem \ref{thm:rcf1},  
to show that two consecutive Right moves from a heap with rcf 1 cannot result in another heap with rcf 1.
 The proof is purely number theoretic and we use a basic result for our sequences that we have not found in the existing literature on Wythoff's sequences. 

\begin{lemma}\label{lem:bbnotabn}
If $B(i) + B(j) = AB(n)$ then $i=j=n=0.$
\end{lemma}

\begin{proof}
It is well-known that all Beatty sequences are sub-additive. More 
precisely, for all $i,j$, 
\begin{align}\label{one}
B(i)+B(j)\in \{B(i+j)-1, B(i+j)\}.
\end{align}
Suppose by way of contradiction that $n>0$ and $ B(i) + B(j) = AB(n)$ (which is also equal to $B(A(n)) + 1$ as shown by Equations (\ref{Ki2}) and (\ref{Ki3})). We divide the proof into two cases:\\

\noindent \emph{Case 1:} If $A(n) \geqslant i+j$, then $B(i) + B(j) = B(A(n)) + 1 \geqslant B(i + j) + 1 > B(i + j)$, a contradiction.\\

\noindent \emph{Case 2:} If $i+j > A(n)$, then either $B(i +j) = B(A(n)) + 2 = B(A(n)+1)$ or $B(i + j) > B(A(n)) + 2$. Supposing the latter, $B(A(n)) + 1 = B(i) + B(j) \geqslant B(i + j) - 1 > B(A(n)) + 1$, a contradiction. Supposing the former, $B(n) = B(n) + A(n) - 1 - A(n) + 1 = B(A(n)) - A(n) + 1 = B(A(n) + 1) - (A(n) + 1) = A(A(n) + 1)$, which contradicts complementarity.
\end{proof}

\subsection{Proofs of partitioning}\label{sec:partitioning}
For an integer $c$ and all nonnegative integers $n$, we denote the sets $A\oplus c = \{A(n) + c : n\geqslant 0\}$ and $B\oplus c = \{B(n) + c : n \geqslant 0\}$ (and analogously for $\ominus$). If $c=1$ we also use the notation $\hat \omega = \omega\oplus 1$, that is $\hat \omega(i)=\omega(i)+1$, for all $i$, if $\omega$ is a word on the alphabet $\{A,B\}$. We also write $\omega_0=\{\omega(i) : i\geqslant 0\}$.

The following lemma proves the first part of Lemma~\ref{lem:part}.
\begin{lemma}[Partitioning Lemma 1]\label{lem:part1} 
The sets $B$, $AB_0$, $\widehat{AB_0}$, and $\widehat{B^2}$ partition the nonnegative integers.
\end{lemma}
\begin{proof}
First subtract 1 from each element of each set; omit $n=0$ in $AB_0$ and $\widehat{AB_0}$. We get the sets $A^2, AB, B^2, BA$, the first by $B\ominus 1=A^2$ and the last by $AB\ominus 1=BA$. 
The former two sets partition $A$, whereas the latter two partition $B$, so the result follows since, shifting back the sets by adding 1, the case $n=0$ provides the ``0" entry from $AB_0$ and the ``1" entry from $\widehat{AB_0}$.
\end{proof}

We prove the second part of Lemma~\ref{lem:part}, restated as Lemma~\ref{lem:part2}, 
via some other partitioning lemmas. 

\begin{lemma}\label{lem:32}
Let $\omega $ be a finite word on the 2-letter alphabet $\{A,B\}$. Then, for $n\geqslant 0$, the sets $\Omega_n=\{\omega B^nA(i) : i>0\}$ partition the set $\Omega=\{\omega(i): i>0\}$.
\end{lemma}
\begin{proof}
Suppose that $x\in \Omega\setminus\cup_{n\geqslant 0}\Omega_n$. By the complementarity of $A$ and $B$, there must exist a smallest $k$ such that $x=\omega B^k(m)$, for some $m$. By minimality of $k$ it follows that $m\not\in B$; hence $m\in A$, which contradicts the definition of $x$. Hence all members in $\Omega$ are represented. Suppose next that there is a smallest $x\in \Omega$ such that $x=\omega B^mA(i)=\omega B^nA(j)$ for some $m< n$ and $i> j$ (since each sequence is increasing we may just as well assume double inequality). By $m<n$ we must then be able to write $A(i)=B^{n-m}A(j)$, which contradicts the complementarity of the $A$ and $B$ sequences.
\end{proof}
We note that if $\omega$ is the empty word, then viewed as a function it is the identity and so, in this case $\Omega$ will be the positive integers. 
For two sets X and Y, we use the notation $X\sqcup Y$ to mean $X\cup Y$, but where we also claim that $X\cap Y=0$, and for several sets the latter equation holds for any pairwise combination.

\begin{lemma}\label{lem:Claim1} $B^3_0\oplus 3\sqcup\{AB^nA(i): n>1, i>0\} = \{G_i(1) : i>0\}.$
\end{lemma}
\begin{proof}
By Theorem \ref{theorem:genK}, for all $i$, $B^{2}(i)+3= \lfloor n\phi\rfloor F_{4}+ iF_{3}+3 = G_i(1)$. Also,
 $B^2\setminus B^3=B^2A$. Hence, to prove Claim 1, it suffices to show that $B^2A\oplus 3=\{AB^nA(i): n>1, i>0\}$. 
Clearly, for all $i$, $AB^2(i) = B^2A(i)+3$. (Since $AB^2(1) = B^2A(1)+3$ and since $AB^2(i)-AB^2(i-1) = B^2A(i) -B^2A(i-1)$ 
for all $i>0$ by (\ref{BnABnA}).) By Lemma \ref{lem:32} we have that $AB^2=\{AB^nA(i): n>1, i>0\}$.
\end{proof}

\begin{lemma}\label{lem} 
$ABA= (B^3_0\oplus3)\sqcup (AB^2\oplus 3).$
\end{lemma}

\begin{proof} We show the equivalent statement $ABA\ominus 3\setminus B^3_0=\{8,21,29,\ldots\}=AB^2$. 
The first equality is clear by Lemma~\ref{lem:Claim1}. We have that $ABA=B^2\ominus 2$. Hence we want to show $B^2\ominus 5\setminus B^3_0=AB^2=B^2A\oplus 3$ which is equivalent to  $B^2\setminus (B^3_0\oplus 5) = B^2A\oplus 8$. Thus, to prove the claim, it suffices to prove 
\begin{align}\label{b3b2a}
B^2 = (B^3_0\oplus 5) \sqcup (B^2A\oplus 8). 
\end{align}
However,
 this will follow from (\ref{BnBn2}). Namely, it suffices to prove $B^3_0\oplus 5 = B^2(B_0\oplus 1)$ and 
 $B^2A\oplus 8 = B^2(A\oplus 1)$, since ($B_0\oplus1)\sqcup (A\oplus 1)=\mathbb Z_{>0}.$ We have that 
 $B^3(0)+ 5 = 5 = B^2(B(0)+1)$ and $B^2A(1)+ 8 =13= B^2(A(1)+1)$ which means that it suffices to prove that, 
 for all $i$, $B^3(i)-B^3(i-1)=B^2(B(i)+1)-B^2(B(i-1)+1)$ and $B^2A(i)-B^2A(i-1)=B^2(A(i)+1)-B^2(A(i-1)+1)$. 
 Therefore (\ref{BnBn2}) and (\ref{BnABnA}) give (\ref{b3b2a}).
\end{proof}

A perhaps less known variation of Lemma \ref{lem:32} is as follows. 
\begin{lemma}\label{lem:32b}
Let $\omega$ denote a finite word on the 2-letter alphabet $\{A,B\}$. Then, for $n\geqslant 0$, the sets $\Omega_n = \{\omega {\hat B}^n\hat A(i) : i\geqslant 0\}$ partition the set $\Omega = \{\omega(i): i>0\}$.
\end{lemma}
\begin{proof}
Note that the sets $\{\hat A(i):i\geqslant 0\}$ and $\{\hat B(i):i>0\}$ are complementary on the positive integers (of course $\hat A(0)=\hat B(0)=1$). Suppose that $x\in \Omega\setminus\cup_{n\geqslant 0}\Omega_n$. Then, by complementarity, there must exist a smallest $k$ such that $x=\omega \hat B^k(m)$, for some $m\ne 1$. By minimality of $k$, $m\not\in \hat B$; hence, again by complementarity, $m\in \hat A$, which contradicts the assumption. The proof of the second part is similar to that of the second part of the proof of Lemma \ref{lem:32}.
\end{proof}

\begin{lemma}\label{lem:33}
For all $i\geqslant 0$ and $n > 0$, $G_{i}(n)=B^{n+1}(i) + F_{2n+3} - 2 = \hat A\hat B^n\hat A(i)-2$.
\end{lemma}
\begin{proof}
Recall, for $i\in \Nz$ and $n\in \N$, $G_i(n):=F_{2n+2}\lfloor i\phi\rfloor + F_{2n+1}i + F_{2n+3} -2.$ 
The first equality is clear by Theorem \ref{theorem:genK}. Since, for a fixed $n$, the gaps of the corresponding sequences are the same, it suffices to check that the first entries of the respective sequences hold, that is that $B^{n+1}(0) + F_{2n+3} - 1=F_{2n+3}-1$ and $ A\hat B^n\hat A(0)=A\hat B^n(1)$ coincide. Since, here, the initial $A$ in the word is not 1-shifted, it suffices to prove that $\hat B^n(1)=F_{2n+2}$, by the third equation in Lemma \ref{lem:aF2n}.
Note $\hat B(1)=2+1=F_4$ and suppose that this holds for some $n>0$. Then $\hat B^{n+1}(1)=\hat B(F_{2n+2})= F_{2n+4}-1+1=F_{2n+4}$, by (\ref{BF2n}). 
\end{proof}

\begin{lemma}\label{lem:new}
$$\bigsqcup_{n>1}G(n)= AB^2\oplus 3.$$ 
\end{lemma}
\begin{proof}
By Lemma \ref{lem:33} we have that, for all $i\geqslant 0$ and (in particular) all $n > 1$, $G_{i}(n)-3 = \hat A \hat B^2\hat B^{n-2}\hat A(i)-5$. Since $B\oplus 1\subset A$, (\ref{BnBn2}) gives that, for all $i$, $\hat A\hat B^2(i)=A(B(B(i)+1)+1)+1=AB^2(i)+5$, which gives the result, by applying $\omega=AB^2$ in Lemma \ref{lem:32b}. 
\end{proof}

\begin{lemma}[Partitioning Lemma 2]\label{lem:part2} 
The sets $G(n)$, for $n>0$, partition the set $AB=\{AB(n): n \geqslant 1\}$.
\end{lemma}
\begin{proof} 

To show that the $G_n$ sets  partition $AB$, we begin by noting that Lemma~\ref{lem} and Lemma~\ref{lem:new} together give 
$(B^3_0\oplus 3)\sqcup\{G_i(n) : n>1, i>0\}=ABA$. Thus, by applying Lemma~\ref{lem:Claim1}, we get 
\begin{align*}
\{G_i(n) : n\geqslant 1, i>0\} &=\{G_i(1) : i>0\}\sqcup \{G_i(n) : n>1, i>0\}\\ 
&= ((B^3_0\oplus 3)\cup\{AB^nA(i): n>1, i>0\}) \sqcup (ABA\setminus (B^3_0\oplus3))\\ 
&= \{AB^nA(i): n\geqslant 1, i>0\}\\
&= AB,
\end{align*}
 since by Lemma \ref{lem:32}, we have that, for $n\geqslant 1$, the sets $AB^nA$ partition $AB$.
 \end{proof}

\section{Discussion}\label{sec:disc}

In Question \ref{qu:example} we asked: Who wins when there are blue heaps of 3 and 20 and a red heap of 18; 
(ii) a blue heap of 20 and a red heap of 17. 

In the first case, the reduced canonical forms are $1/2$, $\{1|0\}$ and $\{0|-1\}$. Since $\{1|0\}+\{0|-1\}=0$, the actual value of a blue
heap of 20 and a red heap of 18 is an infinitesimal. This added to $1/2$ is positive and so Left wins going first (e.g. remove 17 from 18)
and second.

In the second case, the reduced canonical form is $\{1|0\}-1$. Right wins going first (remove the 20 heap). Left going first goes to 0 and
 again the actual value is an infinitesimal but we do not know the sign or whether it is zero. In this case we have to revert to the 
 canonical forms of the games. From Table  \ref{tab:can_red}, we see that Left can remove 16 from 20 to leave a combined value of 
 $\{1||1|0\}-1=1+\{0||0|-1\}-1 = \{0||0|-1\} >0$ and so wins.
 
 Any partizan subtraction game with complementary subtraction sets we call a \cgname{ComplementarySubtraction}, 
CS, game. Any set of positive integers defines a game in CS.
If one of the player's, say Right's, subtraction set is finite, the asymptotic behaviour is obvious, namely, for large heap sizes, Left wins and the value tends to arbitrarily large numbers. Therefore, as in \WPS\!\!, it makes more sense that both players have infinite subtraction sets.

We alluded to the outcome sequence for \wps\  in the Introduction. However, the description applies to any 
CS game based on Beatty sequences.
\begin{theorem}\label{thm:outcomes}
Let $\alpha\in(1,2)$ be irrational and $A(\alpha), B(\alpha)$ be the associated complementary Beatty sequences. In the corresponding
 \cgname{ComplementarySubtraction} game, Left can remove any number in $A(\alpha)$ and Right any number in $B(\alpha)$.
 In this game, the outcome of a heap of size $n$ is: (i) a Left win if  $n\in A(\alpha)$; and (ii) a next player win if $n\in B(\alpha)$.
\end{theorem}
\begin{proof} If  $n\in A$ then Left wins be removing all the heap. If Right moves first then he leaves a number which is either in 
$A(\alpha)$ or in $B(\alpha)$. If the first case holds then again Left removes the whole heap and wins. In the second case, since there are no consecutive numbers in $B(\alpha)$ then Left can reduce the heap down to $1$ in which Right has no move since $1\in A(\alpha)$. 
If $n\in B$ then Right wins immediately. If Left moves then he leaves a number in $A$ and so she wins.

Similarly, if $n\in B(\alpha)$ then Right wins by removing the whole pile. Left, moving first, can reduce the heap to 1 and win.
\end{proof}

\begin{example}\label{ex:1}
Let Left's and Right's subtraction sets be the odd and even positive integers respectively. 
That is, given a heap size $h$, Left can move to $h-o\ge 0$ and right can move to $h-e\ge 0$, where $o$ and
 $e$ is any odd and even positive integer respectively. Then $h=1$ has canonical form the number $1$ (Left wins) 
 and $h=2$ has canonical form $\cgbgame{1}{0}$ (Next player wins). In general one can see that the outcome is 
 Left wins from odd heap sizes and the Next player wins from even heap sizes. (\WPS shares this property in the 
 sense that Left wins moving from a heap size in her subtraction set and the Next player wins moving from a heap size
  in the complement set.) Hence for this game the outcome function is periodic. (For \wps, it follows that the outcome 
  function is aperiodic, but nearly linear in the sense of Wythoff's sequences.) In general, the canonical form for this game 
  is $1/2^{(h-1)/2}$ for odd $h$ and otherwise, for even heap sizes $h\ge 2$, $f$ is defined by the recurrence 
  $f(h)=\cgbgame{1}{0, f(h-2)}$, $f(0)=0$.
\end{example}
Motivated by Example~\ref{ex:1} and \cite{S} we ask the following question.
\begin{question}
Is it true that, for \cgname{ComplementarySubtraction}, the outcome function is periodic if the subtraction sets are? More precisely, 
let $n>1$ be an integer and define a finite set of positive integers $S$ smaller than $n$. Define another game in \cgname{CS} as: Left's legal moves are congruent to $s$ modulo $n$ for any $s\in S$ and Right's legal moves are congruent to $r$ modulo $n$ for any $r\not\in S$. Is the outcome function necessarily periodic?
\end{question}

\begin{problem} 
Classify those games in CS that share the property described in Theorem~\ref{thm:outcomes} and Example~\ref{ex:1}.
\end{problem}

\end{document}